  \documentclass[prl,preprint,onecolumn,preprintnumbers]{revtex4}
  \usepackage{graphicx}
  \usepackage{dcolumn}
  \newcommand{\beq}{\begin{equation}}
  \newcommand{\eeq}{\end{equation}}
  \newcommand{\beqa}{\begin{eqnarray}}
  \newcommand{\eeqa}{\end{eqnarray}}

  \usepackage{amssymb}
  \usepackage{amsfonts}
  \usepackage{amsmath}
  \usepackage{float}
  \usepackage[all]{xy}
  \usepackage{psfrag}
  \usepackage{amsthm}
  \usepackage{amscd}
  \usepackage{tabularx}
  \usepackage{amstext}
  \usepackage{amsxtra}
  \usepackage{plain}  
\begin{document}

\title{Signature-inverse Theorem in Mesh Group-planes}

\author{\vspace{4mm} Reza Aghayan\\ Department of Mathematics - University of Texas at San Antonio - TX 78249 - USA. \\
Email: Reza.Aghayan@utsa.edu. 
\vspace{8mm}\\}

%\affiliation{Faculty of Basic Sciences, Azad University, Tehran, IRAN\\
%Email:r.aghayan@iiau.ac.ir.
%\\}
%-------------------------------------------------------------------------------------------------------------------------------------------------------------

\begin{abstract}
\vspace{3mm}

This is the second paper devoted to the numerical version of Signature-inverse Theorem in terms of the underlying joint invariants. Signature Theorem and its Inverse guarantee any application of differential invariant signature curves to the invariant recognition of visual objects. We first show the invalidity of Curvature-inverse and Signature-inverse theorems, meaning non-congruent meshes may have the same joint invariant numerical curvature or signature. Then by classifying three and five point ordinary meshes respectively in the Euclidean and affine cases, we look for conditions in terms of the associated joint invariant signatures which make these theorems correct. Additionally, we bring forward The Host Theorem to provide a simpler version of Signature-inverse Theorem for closed ordinary meshes.
\vspace{3mm}

Mathematical subject classification 2010: \textit{53A55, 53A04, 53A15, 14L24, 65D18}. 
\vspace{1.5mm}

Keywords: \textit{Invariants theory, Differential invariant signature curves, Joint invariant numerical signatures, Curvature-inverse and Signature-inverse theorems, and Curve analysis}.

\end{abstract}
\maketitle
%-----------------------------------------------------------------------------------------------------------------------------

\section{1 Introduction}

Geometric invariants play a crucial role in object recognition where the object of interest is affected by a transformation group. They were studied by Halphen \cite{hal}, Wilczynski \cite{wil1,wil2}, \v{C}ech and Fubini \cite{fub}, Weyl \cite{wey}, Cartan \cite{car}, Nagata \cite{nag}, and Mumford \cite{mum} who developed the theory of the invariants of transformation groups. A more modern approach has being studied by Calabi et al. \cite{cal3} by introducing the invariant signatures of planar curves, which later modified by Boutin \cite{bou}, Aghayan et al. \cite{agh}, and Aghayan \cite{agh1, agh2, agh3, agh4, agh5}. In \cite{agh, agh1, bru3, bru2, bru4, bou, fen, cal3, hof2}, differential and integral signatures were applied to invariant recognition of object boundaries and detecting symmetries. 

Throughout this paper $\mathrm{\Bbb{G}}$ refers to the special Euclidean motions $\mathrm{S\Bbb{E}}$, the Euclidean group $\mathrm{\Bbb{E}}$, the equiaffine transformations $\mathrm{S\Bbb{A}}$, and the extended equiaffine group $\mathrm{\bar{\Bbb{A}}}$. Also, a $\mathrm{\Bbb{G}}$-plane $\mathrm{E^{\Bbb{G}}}$ means the plane $\mathrm{E\simeq \Bbb{R}^{2}}$ with the geometry induced by $\mathrm{\Bbb{G}}$ acting on E.

Curvature Theorem and its Inverse \cite{sap} indicates that a $\mathrm{C^{r}}$ $\mathrm{(r\geq 2}$ in $\mathrm{E^{\Bbb{E}}}$ and $\mathrm{r\geq 4}$ in $\mathrm{E^{\Bbb{A}})}$ curve $\mathrm{\gamma\subset E^{\Bbb{G}}}$ is uniquely represented, up to the transformation group $\mathrm{\Bbb{G}}$, by its $\mathrm{\Bbb{G}}$-invariant curvature $\mathrm{\kappa_{\Bbb{G}}(\iota)}$ as a function of the $\mathrm{\Bbb{G}}$-invariant arc length $\iota$. To avoid the ambiguity caused by the choice of initial point from where the arc length is measured, Calabi et al. \cite{cal3} brought forward differential invariant signature curves (DISCs) or classifying curves as a new scheme of the invariant recognition for visual objects.

\textbf{Definition 1.1 \cite{agh}.} \textit{A $C^{r}$ $(r\geq 3$ in $\mathrm{E^{\Bbb{E}}}$ and $r\geq 5$ in $\mathrm{E^{\Bbb{A}})}$ curve $\mathrm{\gamma\subset E^{\Bbb{G}}}$ is regular if $\mathrm{\kappa_{\Bbb{G}}(\iota)}$ and its first derivative $\mathrm{\kappa_{\Bbb{G},\iota}}$ are defined and analytic over $\mathrm{\gamma}$. Then,} the $\mathrm{\Bbb{G}}$-invariant signature set \textit{of the regular curve $\mathrm{\gamma}$ is parameterized by}
\vspace{-3.7mm}
\begin{eqnarray*}
\mathrm{\Xi_{\Bbb{G}}(\gamma)=\lbrace (\kappa_{\Bbb{G}}(x),\kappa_{\Bbb{G},\iota}(x)) \mid x\in \gamma \rbrace}\subset \Bbb{R}^{2}.
\end{eqnarray*}

\vspace{-3.7mm}
\noindent \textit{Moreover, if $\mathrm{\gamma}$ is nonsingular, i.e. its signature set is a nondegenerate curve, $\mathrm{\Xi_{\Bbb{G}}(\gamma)}$ is called} the $\mathrm{\Bbb{G}}$-DISC \textit{of} $\mathrm{\gamma}$. \textit{Nonsingularity is guaranteed by $\mathrm{(\kappa_{\Bbb{G},\iota},\kappa_{\Bbb{G},\iota\iota})\neq 0}$}.

\textbf{Signature Theorem in $\mathrm{E^{\Bbb{G}}}$ \cite{olv3}.} \emph{Let $\mathrm{\gamma, \tilde{\gamma}\subset E^{\Bbb{G}}}$ be two congruent curves, i.e. $\mathrm{\tilde{\gamma}=g\cdot \gamma}$ for some $\mathrm{g\in \Bbb{G}}$. Then, their DISCs are identical: $\mathrm{\Xi_{\Bbb{G}}(\tilde{\gamma})=\Xi_{\Bbb{G}}(\gamma)}$}.

\textbf{Signature-inverse Theorem in $\mathrm{E^{\Bbb{G}}}$ \cite{olv3}.} \emph{All smooth nonsingular curves with the same DISC are congruent}. 

Accordingly, DISCs can be applied to program a computer to recognize curves modulo a certain transformation group. However, one major difficulty has been the noise sensitivity of standard differential invariants owing to their dependence on high order derivatives. Aiming to obtain less sensitive approximations, Calabi et al. \cite{cal3} suggested numerical expressions for $\mathrm{\kappa_{G}}$ and $\mathrm{\kappa_{G,\iota}}$ in terms of joint invariants and introduced ``joint invariant numerical signatures" (JINSs). Later, Boutin \cite{bou} corrected and Aghayan et al. \cite{agh} generalized the original formulae and recently Aghayan \cite{agh1} illustrated the resulting formulation depends on the viewpoint and introduced `orientation-invariant' JINSs, leading to the same signature for congruent meshes - named ``the current formulation". 

In the first paper in this series \cite{agh3}, we proved that Signature-inverse theorem is not correct in terms of the current formulation and therefore non-congruent meshes may have the same JINS. To deal with the problem, we introduced ``the new formulation" for JINSs and showed that, compared to the current expressions, the new ones are not only closer to $\mathrm{\Bbb{G}}$-DISCs but for ordinary meshes are also more stable. 

This paper is organized as follows. Section 2 provides a brief survey of the framework of the new formulation. Section 3 shows that Curvature-inverse Theorem and Signature-inverse Theorem are not valid in terms of the new expressions, therefore, non-congruent meshes may have identical numerical curvatures or JINSs. Section 4 first classifies equally and unequally spaced three-point ordinary meshes with respect to their curvatures and side lengths, then, we look for conditions in terms of the new formulas to make the Euclidean Signature-inverse Theorem correct. Next, we bring forward The Host Theorem to give simpler versions of this theorem for closed meshes. Section 5 goes through the same process for the equi-affine case by classifying five-point ordinary meshes. Finally, Section 6 present our conclusions.

%
%
%-----------------------------------------------------------------------------------------------------------------------------------------------------------------------------
% 
\vspace{-4mm}
\section{2 K-point $\mathrm{\Bbb{G}}$-signatures - The new formulation}
\vspace{-2mm}

According to \cite{agh2}, the following subsection does a brief survey of the new formulation.

In a mesh of points $\mathrm{\lbrace p_{i} \rbrace \subset E}$, a \textit{cusp} is a point $\mathrm{p_{i}}$ where the moving point $\mathrm{p_{i+1}}$ starts to move backward - in other word, $\mathrm{p_{i+1}=p_{i-1}}$. An \textit{ordinary} mesh point $\mathrm{\gamma^{\vartriangle}=\lbrace p_{i} \rbrace\subset E}$ refers to a set of successive points with no cusp. An ordinary mesh $\mathrm{\gamma^{\vartriangle}}$ is \textit{fine} if it approximates a  $\mathrm{C^{r}} \hspace{1.2mm} \mathrm{(r\geq3}$ in $\mathrm{E^{\Bbb{E}}}$ and $\mathrm{r\geq5}$ in $\mathrm{E^{\Bbb{A}})}$ curve $\mathrm{\gamma^{\vartriangle}\subset \gamma\subset E}$ and all angles $\mathrm{\theta_{i}=\hspace{1.2mm}<\hspace{-1.3mm}(p_{i-1}, p_{i}, p_{i+1})}$ are obtuse - in this case, the JINS of $\mathrm{\gamma^{\vartriangle}}$ approximates the DISC of $\mathrm{\gamma}$.

\noindent Also, \textit{the n-neighborhood} of $\mathrm{p_{i}\in \gamma^{\vartriangle}}$ means the 2n+1 successive points: 
\vspace{-3mm}
\begin{eqnarray*}
\mathrm{p_{i-n}, \ldots, p_{i-1}, p_{i}, p_{i+1}, \ldots, p_{i+n}}.
\end{eqnarray*}

\vspace{-2mm}
\noindent In addition, a ``k-point invariant" of $\mathrm{E^{\Bbb{G}}}$ is a function $\mathrm{J:E^{k}\longrightarrow\Bbb{R}}$ such that for each $\mathrm{g\in \Bbb{G}}$ and every k-point subset $\mathrm{\lbrace p_{1}, \ldots , p_{k}\rbrace \subset E}$ we have $\mathrm{J(g\cdot p_{1},\ldots, g\cdot p_{k})=J(p_{1}, \ldots, p_{k})}.$
\clearpage
%
%
%----------------------------------------------------------------------------------------------------------------------------------------------------------------
%
%
\subsection{2.1 Two-point $\mathrm{S\Bbb{E}}$-signatures}
\vspace{-3mm}

Consider the Euclidean geometry of curves in the plane $\mathrm{E^{S\Bbb{E}}\simeq \Bbb{R}^{2,S\Bbb{E}}}$  where the underlying group of transformations is the special Euclidean group $\mathrm{S\Bbb{E}(2)=SO(2)\ltimes\Bbb{R}^{2}}$, containing all translations and rotations. One can also include reflections, leading to the Euclidean group $\mathrm{\Bbb{E}(2)=O(2)\ltimes\Bbb{R}^{2}}$.

\textbf{Definition 2.1.} The 2-point $\mathrm{\Bbb{E}}$-curvature \textit{of an ordinary mesh $\mathrm{\gamma^{\vartriangle}=\lbrace p_{i} \rbrace_{i\geq3}\subset E^{\Bbb{E}}}$ is the function $\mathrm{\kappa_{\Bbb{E}}^{\vartriangle}:\gamma^{\vartriangle}\longrightarrow \Bbb{R}}$ given by $\mathrm{p_{i}\longmapsto \frac{4\bigtriangleup}{abc}(p_{i})}$, where $\mathrm{a(p_{i})\geqslant b(p_{i})\geqslant c(p_{i})}$ are the Euclidean distances in the one-neighborhood of $\mathrm{p_{i}}$ and $\mathrm{\bigtriangleup(p_{i})}$ is the area of the triangle whose vertices are the one-neighborhood of $\mathrm{p_{i}}$ given by}
\vspace{-3mm}
\begin{eqnarray*}
\mathrm{\bigtriangleup(p_{i})=\sqrt{(a+(b+c))(c-(a-b))(c+(a-b))(a+(b-c))}(p_{i})}.
\end{eqnarray*}

\vspace{-3mm}
\noindent Then \textit{the 2-point $\mathrm{S\Bbb{E}}$-signature} of an ordinary mesh $\mathrm{\gamma^{\vartriangle}=\lbrace p_{i} \rbrace}$ is parameterized as follows.
 
- Where $\mathrm{\gamma^{\vartriangle}}$ is equally spaced:
\vspace{-3mm}
\begin{eqnarray}
\mathrm{\Xi_{S\Bbb{E}}^{\vartriangle}(\gamma^{\vartriangle})=\lbrace (\frac{4\Delta}{abc}(p_{i}), \space\ \frac{\frac{4\Delta}{abc}(p_{i+1})-\frac{4\Delta}{abc}(p_{i})}{d_{i}})\mid p_{i}\in \gamma^{\vartriangle} \rbrace \subset \Bbb{R}^{2}}
\end{eqnarray}

\vspace{-3mm}
\noindent or, in terms of a centered difference quotient
\vspace{-3mm}
\begin{eqnarray}
\mathrm{\Xi_{S\Bbb{E}}^{\vartriangle}(\gamma^{\vartriangle})=\lbrace (\frac{4\Delta}{abc}(p_{i}), \space\ \frac{\frac{4\Delta}{abc}(p_{i+1})-\frac{4\Delta}{abc}(p_{i-1})}{d_{i-1, i+1}})\mid p_{i}\in \gamma^{\vartriangle} \rbrace \subset \Bbb{R}^{2}}.
\end{eqnarray}

\vspace{-3mm}
 - Where $\mathrm{\gamma^{\vartriangle}}$ is unequally spaced: 
\vspace{-3mm}
\begin{eqnarray}
\mathrm{\Xi_{S\Bbb{E}}^{\vartriangle}(\gamma^{\vartriangle})=\lbrace (\frac{4\bigtriangleup}{abc}(p_{i}), \space\ 3\cdot \frac{\frac{4\bigtriangleup}{abc}(p_{i+1})-\frac{4\bigtriangleup}{abc}(p_{i})}{d_{i-1, i+2}})\mid p_{i}\in \gamma^{\vartriangle} \rbrace \subset \Bbb{R}^{2}}
\end{eqnarray}

\vspace{-3mm}
\noindent or, in terms of a centered difference quotient
\vspace{-3mm}
\begin{eqnarray}
\mathrm{\Xi_{S\Bbb{E}}^{\vartriangle}(\gamma^{\vartriangle})=\lbrace (\frac{4\bigtriangleup}{abc}(p_{i}), \space\ 3\cdot \frac{\frac{4\bigtriangleup}{abc}(p_{i+1})-\frac{4\bigtriangleup}{abc}(p_{i-1})}{d_{i-3, i+3}})\mid p_{i}\in \gamma^{\vartriangle} \rbrace \subset \Bbb{R}^{2}}.
\end{eqnarray}
In which $\mathrm{d_{i,j}=\vert p_{i}-p_{j}\vert}$ denotes their Euclidean distance and $\mathrm{d_{i}=d_{i,i+1}}$.
%
%
%====================================================================================================================
%
%
\vspace{-5mm}
\subsection{2.2 Three-point $\mathrm{S\Bbb{A}}$-signatures}
\vspace{-3mm}

Affine geometry is the study of geometric properties of the objects in the plane $\mathrm{E^{S\Bbb{A}}}$ which remain unchanged by signed area-preserving affine transformations $\mathrm{x\longmapsto Ax+b, A\in SL(2)}$ and $\mathrm{x,b\in \Bbb{R}^2}$ - called the equiaffine group $\mathrm{S\Bbb{A}(2)=SL(2)\ltimes\Bbb{R}^{2}}$. One can include reflections, leading to the extended equiaffine group $\mathrm{\bar{\Bbb{A}}(2)}$.

\noindent \hspace{4.1mm}\textbf{Definition 2.2 \cite{agh2}.} The 3-point $\mathrm{\Bbb{A}}$-curvature \textit{of an ordinary convex (no three are collinear) mesh $\mathrm{\gamma^{\vartriangle}=\lbrace p_{i} \rbrace_{i\geq5}\subset E^{\Bbb{A}}}$ is the real-valued function $\mathrm{\kappa_{\Bbb{A}}^{\vartriangle}:\gamma^{\vartriangle}\longrightarrow \Bbb{R}}$ given by $\mathrm{p_{i}\longmapsto\frac{S}{F^{(2/3)}}(p_{i})}$, where F and S are the first and second affine invariants of $\mathrm{\gamma^{\vartriangle}}$.}

Moreover, the 3-point affine arc lengths of $\mathrm{\gamma^{\vartriangle}=\lbrace p_{i} \rbrace\subset E^{S\Bbb{A}}}$ are obtained as follows.

\textbf{Theorem 2.3 \cite{agh2}.} \textit{In the five-neighborhood of any point $\mathrm{p_{i}\in \gamma^{\vartriangle}}$,} the 3-point $\Bbb{A}$-arc length $\mathrm{L_{k,l}=L(p_{k},p_{l})}$ \textit{from $\mathrm{p_{k}}$ to $\mathrm{p_{l}}$ is computed as follows}.

-\textit{If $\mathrm{\frac{S}{F^{2/3}}(p_{i})\neq 0}$, then}
\vspace{-3mm}
\begin{eqnarray*}
\mathrm{L_{k,l}=\vert \frac{S}{F^{2/3}}(p_{i})\cdot [klO] \vert; \hspace{4mm} \textit{where} \hspace{2mm} O=(\frac{\mathrm{BE-CD}}{AC-B^{2}},-\frac{\mathrm{AE-BD}}{AC-B^{2}})(p_{i}),}
\end{eqnarray*}

\vspace{-3mm}
\noindent \textit{and $\mathrm{[klO]}$ equals the signed area of the parallelogram whose sides are $\mathrm{p_{k}-p_{l}}$ and $\mathrm{p_{k}-O}$}.

-\textit{If $\mathrm{\frac{S}{F^{2/3}}(p_{i})=0}$, then}
\vspace{-3mm}
\begin{eqnarray*}
\mathrm{L_{k,l}=\sqrt[3]{\frac{\mathrm{A}^{2}}{\mathrm{AE-BD}}}\left[(x_{k}-x_{l})+\frac{B}{A}(y_{k}-y_{l})\right](p_{i}); \hspace{3mm} \textit{where} \hspace{2mm} (x_{k},y_{k})=p_{k}, \hspace{1mm} (x_{l},y_{l})=p_{l},}
\end{eqnarray*}

\vspace{-3mm}
\noindent \textit{and A, B, C, D, and E are the affine functions of $\mathrm{\gamma^{\vartriangle}}$. Moreover, $\mathrm{L_{k}=L(p_{k},p_{k+1})}$}.

\noindent Then \textit{the 3-point $\mathrm{S\Bbb{A}}$-signature} of the convex mesh $\mathrm{\gamma^{\vartriangle}=\lbrace p_{i}\rbrace\subset M^{S\Bbb{A}}}$ is obtained as follows.

- Where $\mathrm{\gamma^{\vartriangle}}$ is equally spaced:
\vspace{-3mm} 
\begin{eqnarray}
\mathrm{\Xi_{S\Bbb{A}}^{\Delta}(\gamma^{\vartriangle})=\lbrace (\frac{S}{F^{2/3}}(p_{i}), \space\ \frac{\frac{S}{F^{2/3}}(p_{i+1})-\frac{S}{F^{2/3}}(p_{i})}{L_{i}})\mid p_{i}\in \gamma^{\vartriangle} \rbrace \subset \Bbb{R}^{2}}
\end{eqnarray}

\vspace{-3mm}
\noindent or, in terms of a centered difference quotient
\vspace{-3mm}
\begin{eqnarray}
\mathrm{\Xi_{S\Bbb{A}}^{\Delta}(\gamma^{\vartriangle})=\lbrace (\frac{S}{F^{2/3}}(p_{i}), \space\ \frac{\frac{S}{F^{2/3}}(p_{i+1})-\frac{S}{F^{2/3}}(p_{i-1})}{L_{i-1,i+1}})\mid p_{i}\in \gamma^{\vartriangle} \rbrace \subset \Bbb{R}^{2}}.
\end{eqnarray}

\vspace{-3mm}
- Where $\mathrm{\gamma^{\vartriangle}}$ is unequally spaced:
\vspace{-3mm}
\begin{eqnarray}
\mathrm{\Xi_{S\Bbb{A}}^{\Delta}(\gamma^{\vartriangle})=\lbrace (\frac{S}{F^{2/3}}(p_{i}), \space\ 5\cdot \frac{\frac{S}{F^{2/3}}(p_{i+1})-\frac{S}{F^{2/3}}(p_{i})}{L_{i-2,i+3}})\mid p_{i}\in \gamma^{\vartriangle} \rbrace \subset \Bbb{R}^{2}}
\end{eqnarray}

\vspace{-3mm}
\noindent Or, in terms of a centered difference quotient
\vspace{-3mm}
\begin{eqnarray}
\mathrm{\Xi_{S\Bbb{A}}^{\Delta}(\gamma^{\vartriangle})=\lbrace (\frac{S}{F^{2/3}}(p_{i}), \space\ 5\cdot \frac{\frac{S}{F^{2/3}}(p_{i+1})-\frac{S}{F^{2/3}}(p_{i-1})}{L_{i-5,i+5}})\mid p_{i}\in \gamma^{\vartriangle} \rbrace \subset \Bbb{R}^{2}}.
\vspace{-2cm}
\end{eqnarray}
%
%
%================================================================================================================
%
%
\vspace{-13mm}
\section{3. Counterexamples in the space of ordinary meshes}
\vspace{-2mm}

Consider the space of all planar ordinary meshes affected by $\mathrm{\Bbb{G}}$, called the Mesh $\mathrm{\Bbb{G}}$-plane $\mathrm{M^{\Bbb{G}}}$. Two n-point ordinary meshes $\mathrm{\gamma^{\vartriangle}=\lbrace p_{i} \rbrace}$ and $\mathrm{\tilde\gamma^{\vartriangle}=\lbrace \tilde{p}_{i} \rbrace}$ in $\mathrm{M^{\Bbb{G}}}$ are \textit{congruent} if and only if there exists a $\mathrm{g\in \Bbb{G}}$ and a permutation $\mathrm{\sigma: \lbrace 1, \ldots, n\rbrace \longrightarrow \lbrace 1, \ldots, n\rbrace}$ such that 
\vspace{-3mm}
\begin{eqnarray} 
\mathrm{\tilde{p}_{i} = g\cdot p_{\sigma (i)} \hspace{5mm} \mbox{for all} \hspace{3mm} 1\leq i\leq n}.
\end{eqnarray}

\vspace{-4mm}
\noindent This identity introduces an equivalence relation on $\mathrm{M^{\Bbb{G}}}$, called \textit{the $\mathrm{\Bbb{G}}$-congruent classes}. From now on, without loss of generality, we suppose that $\mathrm{\gamma^{\vartriangle}}$ and $\mathrm{\tilde\gamma^{\vartriangle}}$ are congruent if and only if there exists $\mathrm{g\in \Bbb{G}}$ such that $\mathrm{\tilde{p}_{i} = g\cdot p_{i}}$ for all $\mathrm{1\leq i\leq n}$.  

%
%-----------------------------------------------------------------------------------------------------------------
%
\vspace{-5mm}
\subsection{3.1 Counterexamples for Curvature-inverse Theorem in $\mathrm{M^{\Bbb{G}}}$}
\vspace{-3mm}

According to \cite{agh, agh1}, Curvature Theorem is correct in $\mathrm{M^{\Bbb{G}}}$, meaning congruent meshes have the same k-point $\mathrm{\Bbb{G}}$-curvature $\mathrm{\kappa_{\Bbb{G}}^{\vartriangle}}$. Now, we are going to investigate its inverse. 
\clearpage

\textbf{Example 1.} Let $\mathrm{\gamma_{0-1}^{\vartriangle}=\lbrace p_{1}, p , p_{2}\rbrace}$ and $\mathrm{\gamma_{0-2}^{\vartriangle}=\lbrace q_{1}, q, q_{2}\rbrace}$ be two 3-point ordinary meshes of the unit circle $\mathrm{\gamma_{0}\subset \Bbb{R}^{2}}$ as shown in Figure 1. 
\begin{figure*}[!hbt]
\vspace{-3.5mm}
\includegraphics[angle=0,scale=0.39]{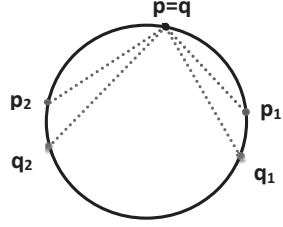}
\vspace{-5.5mm}
\caption{Two arbitrary 3-point ordinary meshes of $\mathrm{\gamma_{0}\subset \Bbb{R}^{2}}$.}
\vspace{-5mm}
\end{figure*} 

\noindent These meshes have the following identical $\mathrm{\Bbb{E}}$-curvatures, while they are clearly non-congruent:

\vspace{-9mm}
\begin{eqnarray*}
\mathrm{\kappa_{\Bbb{E}}^{\vartriangle}(\gamma_{0-1}^{\Delta})=\kappa_{\Bbb{E}}^{\vartriangle}(\gamma_{0-2}^{\Delta}): \lbrace p=q\rbrace \longmapsto \lbrace 1\rbrace}.
\end{eqnarray*}

\vspace{-2mm}
\noindent \hspace{3.4mm} \textbf{Example 2.} Let $\mathrm{\gamma_{1-1}^{\vartriangle}=\lbrace p_{1}, p, p_{2}\rbrace\subset \gamma_{1}}$ and $\mathrm{\gamma_{2-1}^{\vartriangle}=\lbrace q_{1}, q, q_{2}\rbrace\subset \gamma_{2}}$ be two equally spaced 3-point ordinary meshes in $\mathrm{E^{S\Bbb{E}}}$, see FIG. 2, and their circumcircles have the same radius $\mathsf{R}$. %Let $\mathrm{d(p, p_{1})=d(q, q_{2})=d}$. 
\vspace{-3mm}
\begin{figure*}[!hbt]
\vspace{-5mm}
\includegraphics[angle=0,scale=0.415]{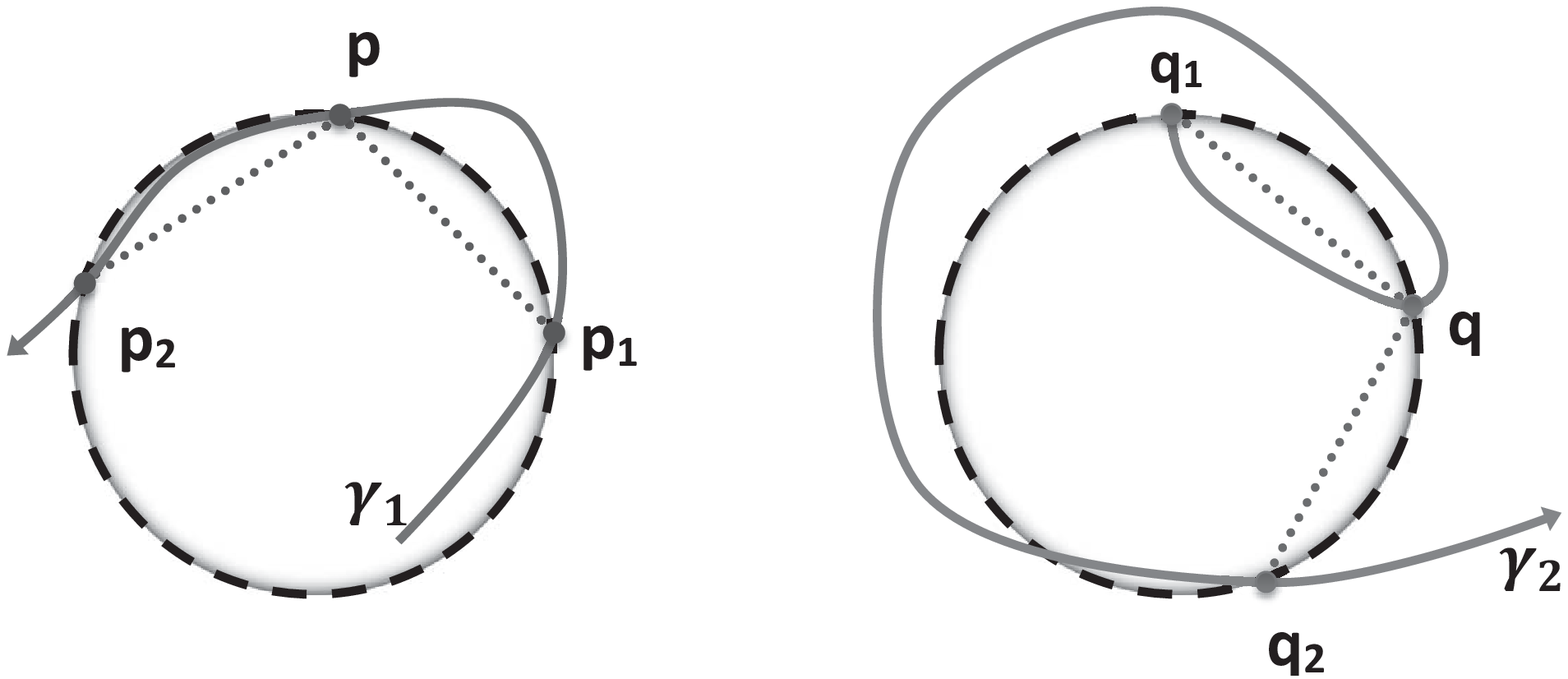}
\vspace{-7mm}
\caption{To transform $\mathrm{\gamma_{2-1}^{\vartriangle}}$ to $\mathrm{\gamma_{1-1}^{\vartriangle}}$ a reflection is required.}
\vspace{-5mm}
\end{figure*} 

\noindent $\mathrm{\gamma_{1-1}^{\vartriangle}}$ and $\mathrm{\gamma_{2-1}^{\vartriangle}}$ are clearly non-congruent, while they have the same 2-point $\mathrm{\Bbb{E}}$-curvature:
\vspace{-3mm} 
\begin{eqnarray*}
\mathrm{\kappa_{\Bbb{E}}^{\vartriangle}(p)= \kappa_{\Bbb{E}}^{\vartriangle}(q) = 1/\mathsf{R}}. 
\end{eqnarray*} 

\vspace{-3mm}
Now consider two approximating meshes of $\mathrm{\gamma_{0}}$ with the same resolution (number of points) in which one is equally and the other one is unequally spaced. They have the same 2-point $\mathrm{\Bbb{E}}$-curvatures, while they are not congruent, meaning Curvature-inverse Theorem is not correct even for \textit{fine} partitions.

Just like the Euclidean case, it is also easy to give some counterexamples in $\mathrm{M^{\Bbb{A}}}$, therefore, we have the following proposition.

\textbf{proposition 3.1.} \textit{Non-congruent meshes may have identical k-point $\mathrm{\Bbb{G}}$-curvatures}.
%
%
%=================================================================================================================
%
%
\vspace{-5mm}
\subsection{3.2 Counterexample for Signature-inverse Theorem in $\mathrm{M^{\Bbb{G}}}$}
\vspace{-4mm}

According to \cite{agh1}, to validate Signature Theorem in terms of the current formulation, one just needs to consider the orientation-invariant version of it, meaning congruent meshes have the same orientation-invariant JINS. Now, we investigate the correctness of its inverse. 

Let $\mathrm{R_{p}}$ be the radius of the circumcircle passing through the one-neighborhood of a point $\mathrm{p}$ on an ordinary mesh $\mathrm{\gamma^{\vartriangle}}$.
\clearpage

\noindent \hspace{4.2mm}\textbf{Example 3.} Let $\mathrm{\gamma_{3-1}^{\vartriangle}=\lbrace p_{1}, p_{2} , p_{3}, p_{4}, p_{5}\rbrace\subset \gamma_{3}}$ and $\mathrm{\gamma_{4-1}^{\vartriangle}=\lbrace q_{1}, q_{2}, q_{3}, q_{4}, q_{5}\rbrace\subset \gamma_{4}}$ denote two equally spaced ordinary meshes with $\mathrm{d(p_{1}, p_{3})=d(q_{1}, q_{3})}$, see FIG. 3. Also, let 

\hspace{2.5cm}$\mathrm{R_{p_{1}}=R_{p_{2}}=}$ $\mathrm{R_{q_{1}}=R_{q_{2}}=\mathsf{R}}$ \hspace{5mm} and \hspace{6mm}$\mathrm{R_{p_{3}}=R_{q_{3}}=\mathsf{r}}$. 
\begin{figure*}[!hbt]
\vspace{-4.3mm}
\includegraphics[angle=0,scale=0.6]{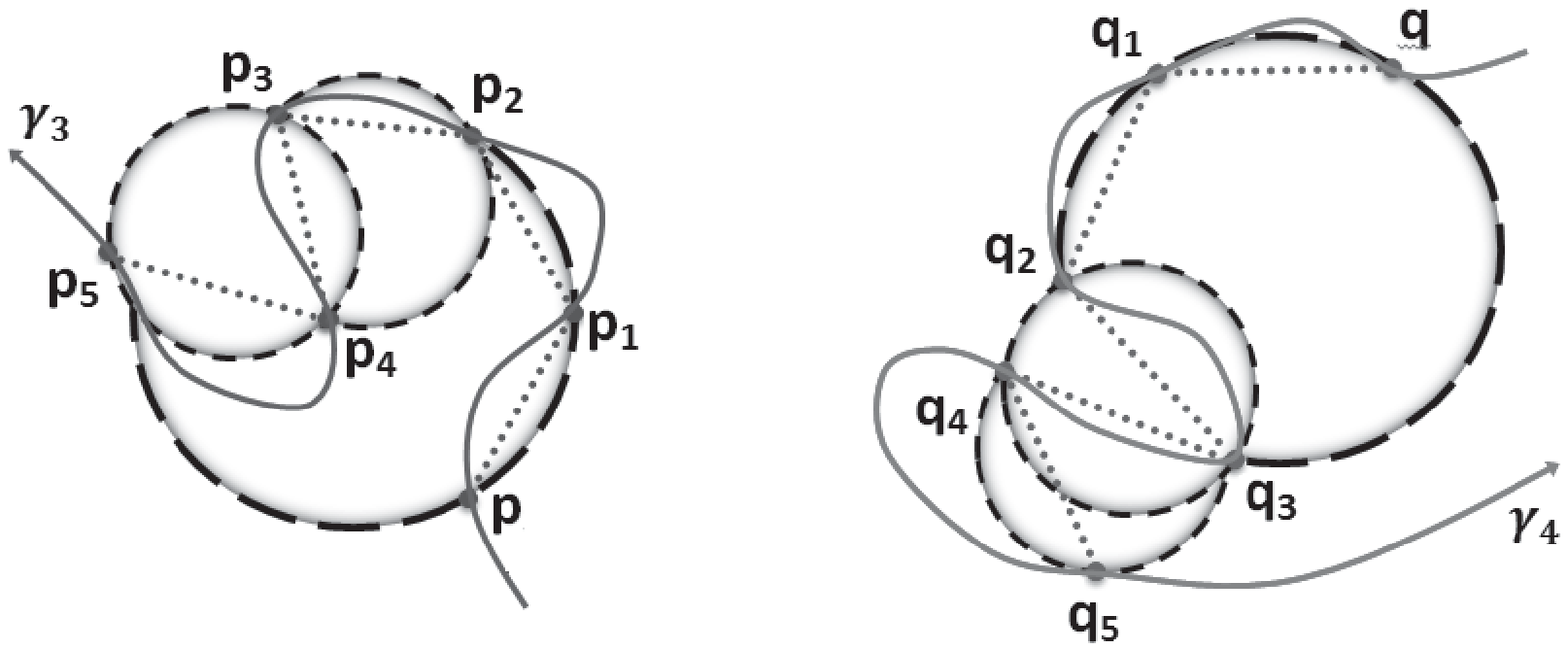}
\vspace{-7mm}
\caption{$\mathrm{\gamma_{3-1}^{\vartriangle}}$ and $\mathrm{\gamma_{4-1}^{\vartriangle}}$.}
\vspace{-7mm}
\end{figure*}

\noindent From identity (2)
\vspace{-3mm}
\begin{eqnarray*}
\mathrm{\Xi_{S\Bbb{E}}^{\vartriangle}(\gamma_{3-1}^{\vartriangle})} &=& \mathrm{\lbrace (\frac{4\Delta}{abc}(p_{i}), \frac{\frac{4\Delta}{abc}(p_{i+1})-\frac{4\Delta}{abc}(p_{i-1})}{d_{i-1, i+1}})\mid p_{i}\in \gamma_{3-1}^{\vartriangle} \rbrace} = \mathrm{\lbrace (\dfrac{1}{\mathsf{R}}, \frac{\frac{1}{\mathsf{r}}-\frac{1}{\mathsf{R}}}{d_{1, 3}})\rbrace}
\end{eqnarray*}

\vspace{-4mm}
\noindent and
\vspace{-4mm}
\begin{eqnarray*}
\mathrm{\Xi_{S\Bbb{E}}^{\vartriangle}(\gamma_{4-1}^{\vartriangle})} &=& \mathrm{\lbrace (\frac{4\Delta}{abc}(q_{i}), \frac{\frac{4\Delta}{abc}(q_{i+1})-\frac{4\Delta}{abc}(q_{i-1})}{d_{i-1, i+1}})\mid q_{i}\in \gamma_{4-1}^{\vartriangle} \rbrace} =\mathrm{\lbrace (\dfrac{1}{\mathsf{R}}, \frac{\frac{1}{\mathsf{r}}-\frac{1}{\mathsf{R}}}{d_{1, 3}})\rbrace}.
\end{eqnarray*}

\vspace{-1mm}
\noindent Thus $\mathrm{\gamma_{3-1}^{\vartriangle}}$ and $\mathrm{\gamma_{4-1}^{\vartriangle}}$ have the same 2-point $\mathrm{S\Bbb{E}}$-signature, while they are not clearly congruent.  

Just like the Euclidean case, it is also easy to give some counterexamples in $\mathrm{M^{S\Bbb{A}}}$, therefore, we have the following corollary.

\textbf{Proposition 3.2.} \textit{Signature-inverse Theorem is not correct in $\mathrm{M^{\Bbb{G}}}$}.

\textbf{Remark.} There are theorems with very simple conditions for Signature-inverse Theorem in $\mathrm{M^{\Bbb{G}}}$, but they are not in terms of the resulting JINSs. For example, the following theorem is the simplest.

\noindent \hspace{2.6mm}\textbf{Theorem 3.3.} \textit{Let $\mathrm{\gamma^{\vartriangle}, \tilde{\gamma}^{\vartriangle}\subset M^{S\Bbb{E}}}$ denote ordinary meshes with the same} Euclidean distance \textit{between any two corresponding successive points and have also identical corresponding} signed angles. \textit{Then, $\mathrm{\gamma^{\vartriangle}}$ and $\mathrm{\tilde{\gamma}^{\vartriangle}}$ are congruent}.

From now on, this paper looks for conditions in terms of the resulting signatures to make Signature-inverse Theorem correct in $\mathrm{M^{\Bbb{G}}}$.
%
%
%=================================================================================================================
%
%
\vspace{-4mm}
\section{4 Signature-inverse Theorem in $\mathrm{M^{S\Bbb{E}}}$}
\vspace{-2mm}

\textbf{Definition 4.1.} \textit{Let $\mathrm{\gamma^{\vartriangle}=\lbrace p_{i} \rbrace\subset E\simeq \Bbb{R}^2}$ denote an ordinary mesh point}. 

\noindent a) The signature-sign \textit{of $\mathrm{\gamma^{\vartriangle}}$ is the function $\mathrm{SS:\gamma^{\vartriangle}\longrightarrow \lbrace-1, 0, 1\rbrace}$ given at any point $\mathrm{p_{i}\in \gamma^{\vartriangle}}$ by}: $\mathrm{p_{i}\longmapsto sgn((p_{i+1}-p_{i})\times(p_{i-1}-p_{i}))}$ \textit{where $\times$ denotes the 2d-cross product}.

\noindent b) \textit{$\mathrm{p_{i}\in \gamma^{\vartriangle}}$ is} in the signature-direction $\mathrm{(SD)}$ \textit{if its one-neighborhood are in counterclockwise order on their circumcircle. Otherwise, we say $\mathrm{p_{i}}$ is in $\mathrm{\sim SD}$}. 

\noindent \textit{Besides, two ordinary meshes $\mathrm{\gamma^{\vartriangle}=\lbrace p_{i} \rbrace_{1}^{n}}$ and $\mathrm{\tilde\gamma^{\vartriangle}=\lbrace \tilde{p}_{i} \rbrace_{1}^{n}}$ are in the same signature-direction if and only if for $1<i<n$ the corresponding points $\mathrm{p_{i}}$ and $\mathrm{\tilde{p}_{i}}$ are both in $\mathrm{SD}$ or in $\mathrm{\sim SD}$}.

\textbf{Proposition 4.2.} \textit{Rotations do not change the signature-directions while reflections do. In other words, for a mesh point $\mathrm{\gamma^{\vartriangle}=\lbrace p_{i} \rbrace}$}

\noindent a) $\mathrm{p_{i}}$ \textit{is in} SD \textit{iff} $\mathrm{\theta\cdot p_{i-1}\curvearrowright \theta\cdot p_{i}\curvearrowright \theta\cdot p_{i+1}}$ \textit{are in} SD, \textit{where} $\mathrm{\theta\in SO(2)}$,

\noindent b) $\mathrm{p_{i}}$ \textit{is in} SD \textit{iff} $\mathrm{\vartheta\cdot p_{i-1}\curvearrowright \vartheta\cdot p_{i}\curvearrowright \vartheta\cdot p_{i+1}}$ \textit{are in} $\mathrm{\sim SD}$, \textit{where} $\mathrm{\vartheta\in}$ $\mathrm{O(2)\diagdown SO(2)}$, 

\noindent meaning the signature-direction makes a distinction between rotations and reflections.
%
%
%=================================================================================================================
%
%
\vspace{-5mm}
\subsection{4.1 Equally spaced meshes in $\mathrm{E^{S\Bbb{E}}}$}
\vspace{-4mm}

A mesh point $\mathrm{\gamma^{\vartriangle}\subset E^{\Bbb{E}}}$ is equally spaced, if all edges (two successive points) have the same Euclidean length. Also, in a 3-point equally spaced mesh, let $\mathrm{\mathsf{d}}$ and $\mathrm{\mathsf{d}\_}$ denote respectively the $\mathrm{\Bbb{E}}$-length of each edge and the $\mathrm{\Bbb{E}}$-distances between the first and end points. In addition, for a n-point equally spaced mesh, let $\mathrm{\mathsf{d}_{i}}$ and $\mathrm{\mathsf{d}_{\_i}}$ be respectively $\mathrm{\mathsf{d}}$ and $\mathrm{\mathsf{d}\_}$ in the one-neighborhood of its $\mathrm{i^{th}}$ point (in open meshes $\mathrm{1<i<n}$).

The following lemma plays a central role in our investigations. 

\textbf{Lemma 4.3 \cite{fel2}.}  \textit{The 2-point action $\mathrm{S\Bbb{E}(2)\times M^{2}\longrightarrow M^{2}}$ given by}\\ 
$\mathrm{\hspace{5.5cm} g\cdot (p_{1}, p_{2})\longmapsto (g\cdot p_{1}, g\cdot p_{2})}$

\noindent \textit{is free on $\mathrm{M^{2}\setminus \mathsf{D}}$, where $\mathrm{\mathsf{D} = \lbrace (p_{1}, p_{2})\mid p_{1}= p_{2}\rbrace}$}.
\vspace{2.5mm}
%
%=================================================================================================================
%

\textbf{$\square \succ$ \textit{Signature-inverse Theorem in terms of (1)}}

\textbf{Lemma 4.4.} \textit{For any $\mathrm{p, q\in E^{S\Bbb{E}}}$, there are at most 2 circles with a given radius $\mathrm{\mathsf{R}}$ passing through them}, see FIG. 4.
\begin{figure*}[!hbt]
\vspace{-10mm}
\includegraphics[angle=0,scale=0.49]{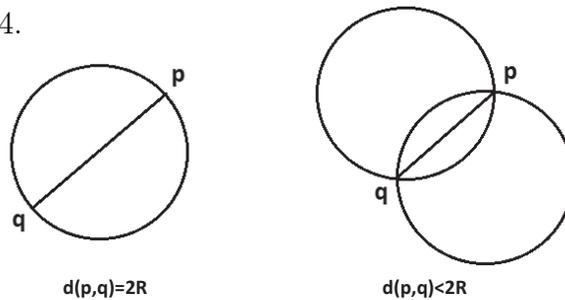}
\vspace{-6.5mm}
\caption{There is no motion in $\mathrm{S\Bbb{E}}$ to coincide the right circles.}
\vspace{-6mm}
\end{figure*} 

\noindent \hspace{2.6mm}Let $\mathrm{\gamma^{\vartriangle}=\lbrace p_{1}, p, p_{2}\rbrace}$ denote a 3-point equally spaced ordinary mesh with $\mathrm{\Bbb{E}}$-curvature $\mathsf{\kappa}$ and side length $\mathsf{d}$, see FIG. 5. 
\begin{figure*}[!hbt]
\vspace{-10mm}
\hspace{-3mm}\includegraphics[angle=0,scale=0.41]{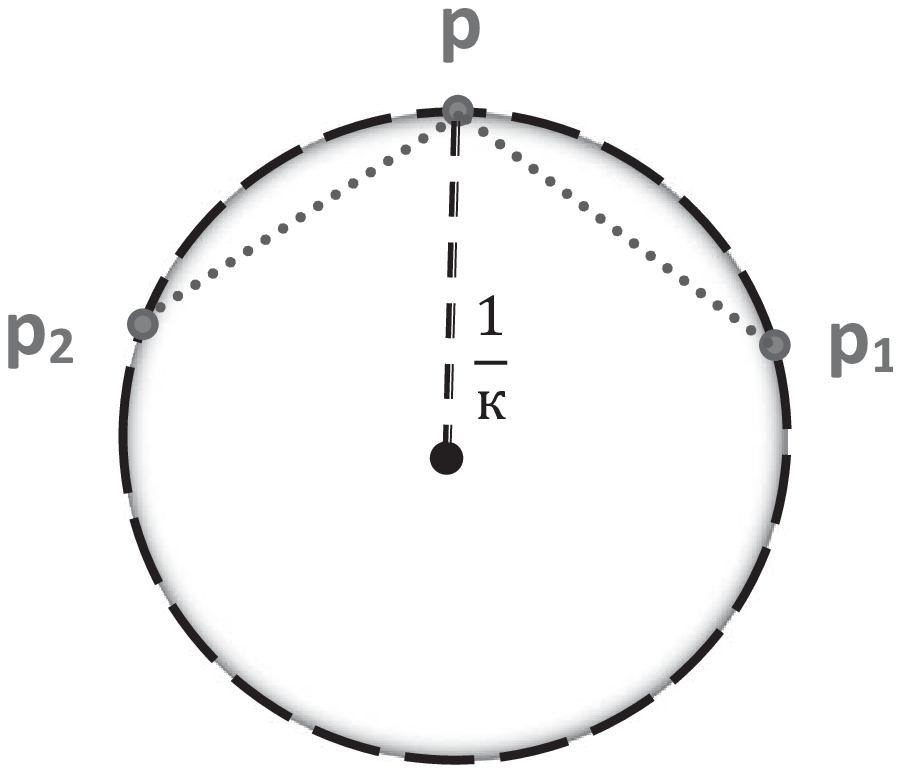}
\vspace{-6.5mm}
\caption{The mesh $\gamma^{\vartriangle}$.}
\vspace{-3.3cm}
\end{figure*} 
\clearpage

\noindent By Lemma 4.4, there exists only one class of 3-point equally spaced meshes with the same $\mathsf{\kappa}$ and $\mathrm{\mathsf{d}}$ and non-congruent with $\gamma^{\vartriangle}$, which is its mirror reflection - denoted by $\mathrm{\tilde{\gamma}^{\vartriangle}=\lbrace q_{1}, q, q_{2}\rbrace}$, see FIG. 6. 
\begin{figure*}[!hbt]
\vspace{-9mm}
\includegraphics[angle=0,scale=0.6]{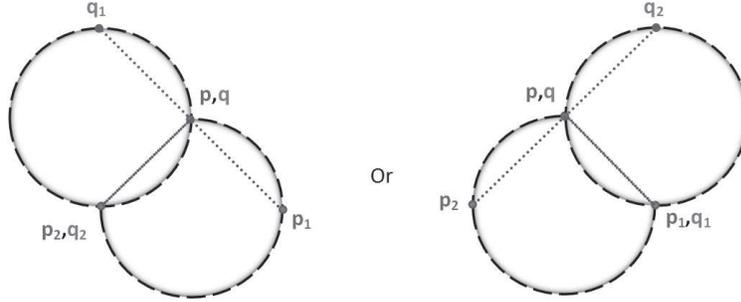}
\vspace{-5mm}
\caption{There is no motion in $\mathrm{S\Bbb{E}}$ to map one mesh to the other.}
\vspace{-6mm}
\end{figure*} 

\noindent \hspace{4mm}We, therefore, have the following results.

\noindent \hspace{4mm}\textbf{Theorem 4.5.} \textit{In $\mathrm{M^{S\Bbb{E}}}$, there exist two congruent classes of equally spaced 3-point ordinary meshes with the same curvature $\kappa$ and side length $\mathsf{d}$. One in $\mathrm{SD}$ and another in $\mathrm{\sim SD}$ denoted respectively by  $\mathrm{\gamma^{\vartriangle}_{\kappa, d, SD}}$ and $\mathrm{\gamma^{\vartriangle}_{\kappa, d, \sim SD}}$, see FIG. 7. Also, there is a unique class of them in $\mathrm{M^{\Bbb{E}}}$}.
\begin{figure*}[!hbt]
\vspace{-10mm}
\includegraphics[angle=0,scale=0.39]{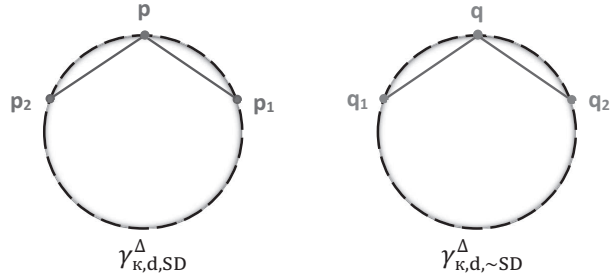} 
\vspace{-6mm}
\caption{Just a reflection in $\mathrm{M^{\Bbb{E}}}$ maps one to the other.}
\vspace{-5mm}
\end{figure*} 

\noindent \hspace{3.1mm}\textbf{Corollary 4.6.} \textit{There exists a unique $\mathrm{S\Bbb{E}}$-congruent class of 3-point equally spaced ordinary meshes with total length $\mathrm{2\mathsf{d}}$, $\mathrm{\Bbb{E}}$-curvature $\kappa$, and being in the same signature-direction}.

\noindent \hspace{2.7mm}\textbf{Corollary 4.7.} \textit{Two equally spaced ordinary meshes $\mathrm{\gamma^{\vartriangle}}$ and $\mathrm{\tilde{\gamma}^{\vartriangle}}$ in $\mathrm{M^{S\Bbb{E}}}$ with the same total length, signature-direction, and 2-point curvature} $\mathrm{\kappa_{\Bbb{E}}^{\vartriangle} ({\tilde\gamma^{\vartriangle}})=\kappa_{\Bbb{E}}^{\vartriangle}({\gamma^{\vartriangle}})}$ \textit{are congruent}.

However, it is \textit{not} in terms of the $\mathrm{S\Bbb{E}}$-JINSs. It is also easy to check the following result.

\textbf{Theorem 4.8.} \textit{In $\mathrm{M^{S\Bbb{E}}}$, there exist two congruent classes of equally spaced 3-point meshes with curvature $\kappa$, angle $0<\theta_{0}<\pi$, and the same signature-direction denoted by $\mathrm{\gamma^{\vartriangle}_{\kappa, \theta_{0}, SD}}$ and $\mathrm{\gamma^{\vartriangle}_{\kappa, \theta_{0}, \sim SD}}$, see FIG. 8. Moreover, there is a unique equivalence class of them in $\mathrm{M^{\Bbb{E}}}$}.
\begin{figure*}[!hbt]
\vspace{-5.7mm}
\includegraphics[angle=0,scale=0.39]{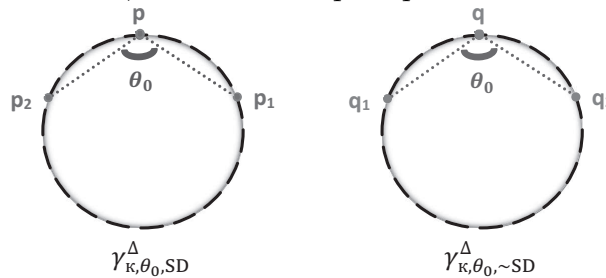} 
\vspace{-5.5mm}
\caption{There is no motion in $\mathrm{S\Bbb{E}}$ to coincide one mesh to the other.}
\vspace{-3.2cm}
\end{figure*}  
\clearpage
%
%--------------------------------------------------------------------------------------------------------------------------
%

\textbf{Theorem 4.9, Signature-inverse Theorem.} \textit{Equally spaced ordinary meshes $\mathrm{\gamma^{\vartriangle}}$ and $\mathrm{\tilde{\gamma}^{\vartriangle}}$ in $\mathrm{M^{S\Bbb{E}}}$ with the same signature-direction and identical JINSs given by (1) are congruent}.

\textbf{Proof.} Let $\mathrm{\gamma^{\vartriangle}=\lbrace p_{i}\rbrace_{1}^{n}}$ and $\mathrm{\tilde{\gamma}^{\vartriangle}=\lbrace \tilde{p}_{i}\rbrace_{1}^{n}}$. Also, without loss of generality, we assume that $\mathrm{\Xi_{S\Bbb{E}}^{\vartriangle}(\tilde{p}_{i})=\Xi_{S\Bbb{E}}^{\vartriangle}(p_{i})}$ for $\mathrm{1\leq i\leq n}$. By identity (1)
\vspace{-3mm}
\[\left\{\begin{array}{cl}
\mathrm{\frac{4\Delta}{abc}(\tilde{p}_{i})=\frac{4\Delta}{abc}(p_{i})}\hspace{3.5cm} \mbox{and}\\ [1.5mm]
\mathrm{\displaystyle \frac{\frac{4\Delta}{abc}(\tilde{p}_{i+1})-\frac{4\Delta}{abc}(\tilde{p}_{i})}{\tilde{d}_{i}}=\frac{\frac{4\Delta}{abc}(p_{i+1})-\frac{4\Delta}{abc}(p_{i})}{d_{i}}}\end{array}\right.\]
%\begin{eqnarray}
%\mathrm{\frac{4\Delta}{abc}(\tilde{p}_{i})=\frac{4\Delta}{abc}(p_{i})} \hspace{4mm} \mbox{and} \hspace{4mm}
%\mathrm{\displaystyle \frac{\frac{4\Delta}{abc}(\tilde{p}_{i+1})-\frac{4\Delta}{abc}(\tilde{p}_{i})}{\tilde{d}_{i}}=\frac{\frac{4\Delta}{abc}(p_{i+1})-\frac{4\Delta}{abc}(p_{i})}{d_{i}}; \hspace{5mm} 1\leq i<n}
%\end{eqnarray}

\vspace{-4mm}
\noindent which result
\vspace{-6mm}
\begin{eqnarray*}
\mathrm{\frac{4\Delta}{abc}(\tilde{p}_{i})=\frac{4\Delta}{abc}(p_{i})=\kappa_{i} \hspace{5mm} \mbox{and} \hspace{5mm} \tilde{d}_{i}=d_{i}; \hspace{5mm} 1\leq i<n}.
\end{eqnarray*}

\vspace{-6mm}
\noindent Hence
\vspace{-6mm}
\begin{eqnarray}
\mathrm{\exists \hspace{.5mm} g_{i}\in S\Bbb{E}(2) \hspace{4mm} \mbox{s.t.} \hspace{4mm} g_{i}\cdot [p_{i}, p_{i+1}]\longmapsto [\tilde{p}_{i}, \tilde{p}_{i+1}]; \hspace{5mm} 1\leq i<n}.
\vspace{-3mm}
\end{eqnarray}

\vspace{-3mm}
\noindent where $\mathrm{[p, q]=tp+(1-t)q; \space\ 0\leq t\leq 1}.$ To show that $\mathrm{g_{i}=g_{i+1}}$ for all $\mathrm{1\leq i<n-1}$, consider the one-neighborhoods of a point $\mathrm{p_{i}}$. From (10)
\vspace{-3mm}
\begin{eqnarray}
\mathrm{g_{i}\cdot [p_{i}, p_{i+1}]\longmapsto [\tilde{p}_{i}, \tilde{p}_{i+1}]} \hspace{4mm} \mbox{and} \hspace{4mm} \mathrm{g_{i+1}\cdot [p_{i+1}, p_{i+2}]\longmapsto [\tilde{p}_{i+1}, \tilde{p}_{i+2}]}.
\end{eqnarray}

\vspace{-4mm}
\noindent On the other hand, Corollary 4.6 and Lemma 4.3 give
\vspace{-4mm}
\begin{eqnarray*}
\mathrm{\exists ! \hspace{.7mm} g^{1}_{i}\in S\Bbb{E}(2) \hspace{4mm} \mbox{s.t.} \hspace{4mm} g^{1}_{i}: N^{1}_{i}\longmapsto \tilde{N}^{1}_{i}; \hspace{5mm} 1<i<n}
\end{eqnarray*} 

\vspace{-4mm}
\noindent which, along with (11), results in 
\vspace{-4mm}
\begin{eqnarray*}
\mathrm{\exists ! \hspace{.5mm} g\in S\Bbb{E}(2) \hspace{4mm} \mbox{s.t.} \hspace{5mm} \tilde{\gamma}^{\Delta}=g\cdot \gamma^{\vartriangle}}.
\end{eqnarray*} 

\vspace{-4mm}
\noindent In other words, $\mathrm{\gamma^{\vartriangle}}$ and $\mathrm{\tilde{\gamma}^{\vartriangle}}$ are $\mathrm{S\Bbb{E}}$-congruent. 

It is not also hard to check that Theorem 4.9 and Theorem 3.3 are the same.
%
%---------------------------------------------------------------------------------------------------------------------------
%

\vspace{2.5mm}
\textbf{$\square \succ$ \textit{Signature-inverse Theorem in terms of (2) - version 1}}
\vspace{2mm}

\textbf{Definition 4.10.} \textit{Two ordinary meshes $\mathrm{\gamma^{\vartriangle}, \tilde{\gamma}^{\vartriangle}\subset {M}}$ have the same ``angle-type" iff any pair of angles $\mathrm{\theta_{i}=\hspace{1.2mm}<\hspace{-1.3mm}(p_{i-1}, p_{i}, p_{i+1})}$ and $\mathrm{\tilde{\theta}_{i}=\hspace{1.2mm}<\hspace{-1.3mm}(\tilde{p}_{i-1}, \tilde{p}_{i}, \tilde{p}_{i+1})}$ are both acute, obtuse, or right}.

\textbf{Theorem 4.11.} \textit{In $\mathrm{M^{S\Bbb{E}}}$, there exist at most 4 congruent classes of 3-point equally spaced ordinary meshes with curvature $\mathrm{\kappa}$ and distance $\mathrm{\mathsf{d\_}}$ as follows}:

- \textit{If} $\mathrm{\kappa\cdot \mathsf{d\_}=2}$; \textit{there are 2 congruent classes}: \textit{one in $\mathrm{SD}$ denoted by $\mathrm{\gamma^{\vartriangle}_{\kappa, SD}}$ and another in $\mathrm{\sim SD}$ denoted by $\mathrm{\gamma^{\vartriangle}_{\kappa, \sim SD}}$, see FIG. 9}.
\begin{figure*}[!hbt]
\vspace{-4mm}
\includegraphics[angle=0,scale=0.39]{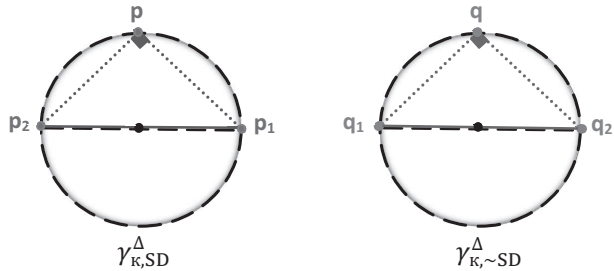}
\vspace{-6mm}
\caption{2 classes of 3-point equally spaced ordinary meshes where $\mathrm{\mathsf{d\_}}$ equals the diameter.}
\vspace{-3.2cm}
\end{figure*} 
\clearpage

- \textit{If} $\mathrm{\kappa\cdot \mathsf{d\_}\neq 2}$; \textit{there are 4 congruent classes}: \textit{two in $\mathrm{SD}$ denoted by $\mathrm{\lbrace\gamma^{\vartriangle,1}_{\kappa, \mathsf{d\_}, SD}, \gamma^{\vartriangle,2}_{\kappa, \mathsf{d\_}, SD}\rbrace}$, and the other two in $\mathrm{\sim SD}$ denoted by $\mathrm{\lbrace\gamma^{\vartriangle,1}_{\kappa, \mathsf{d\_}, \sim SD}, \gamma^{\vartriangle,2}_{\kappa, \mathsf{d\_} \sim SD}\rbrace}$, see FIG. 10}.
\begin{figure*}[!hbt]
\vspace{-4mm}
\includegraphics[angle=0,scale=0.39]{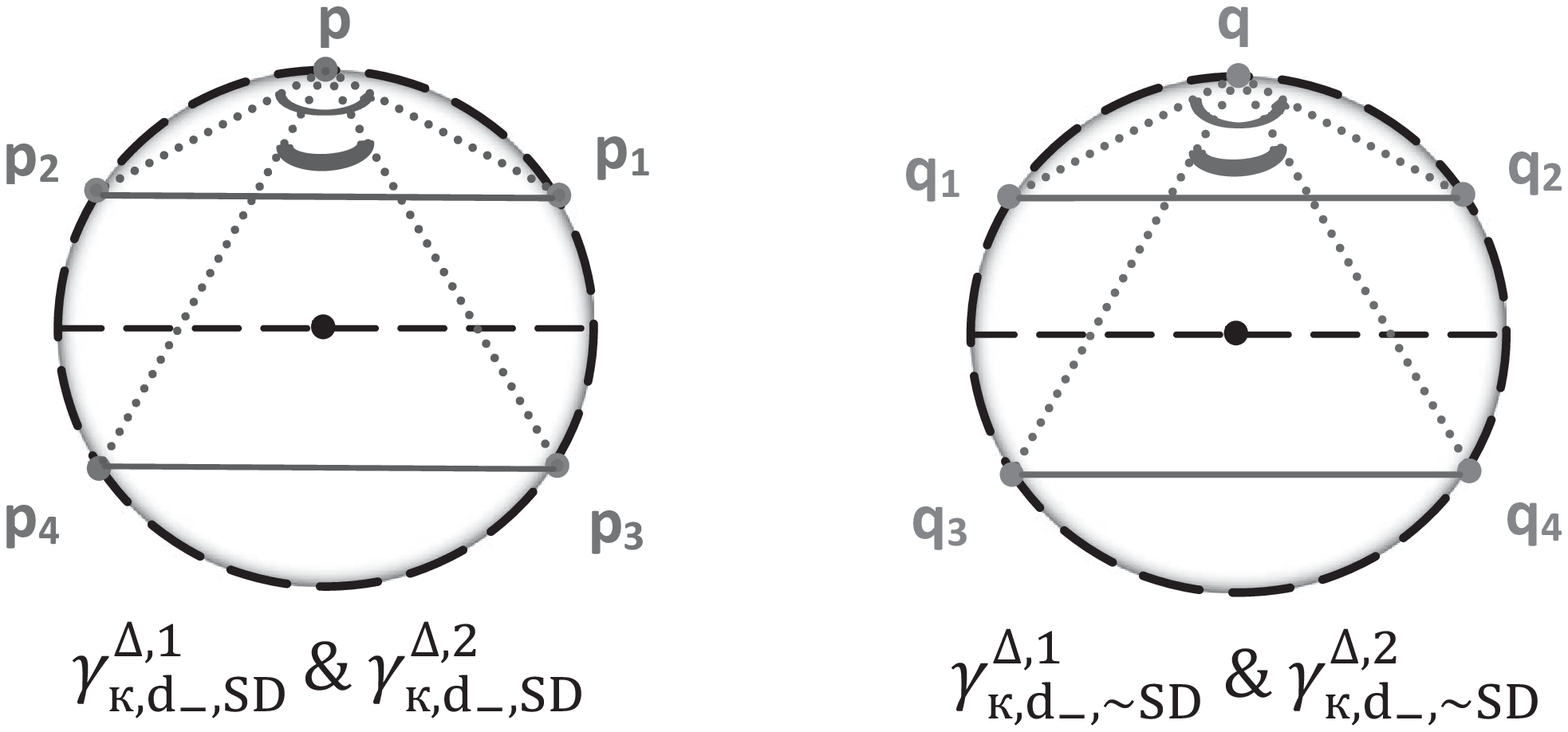}
\vspace{-6mm}
\caption{4 classes of 3-point equally spaced ordinary meshes where $\mathrm{\mathsf{d\_}}$ is less than the diameter.}
\vspace{-4mm}
\end{figure*}  

We, therefore, established the following results.

\textbf{Theorem 4.12.} \textit{In $\mathrm{M^{S\Bbb{E}}}$, there exist 2 congruent classes of 3-point equally spaced ordinary meshes with curvature $\mathrm{\kappa}$, distance $\mathrm{\mathsf{d\_}}$, and also the same angle-type $\mathrm{\theta}$. One in $\mathrm{SD}$ and another in $\mathrm{\sim SD}$ denoted respectively by $\mathrm{\gamma^{\vartriangle}_{\kappa, \theta, SD}}$ and $\mathrm{\gamma^{\vartriangle}_{\kappa, \theta, \sim SD}}$, see FIG. 11. Moreover, there is a unique congruent class of them in $\mathrm{M^{\Bbb{E}}}$}.
\begin{figure*}[!hbt]
\vspace{-5mm}
\includegraphics[angle=0,scale=0.39]{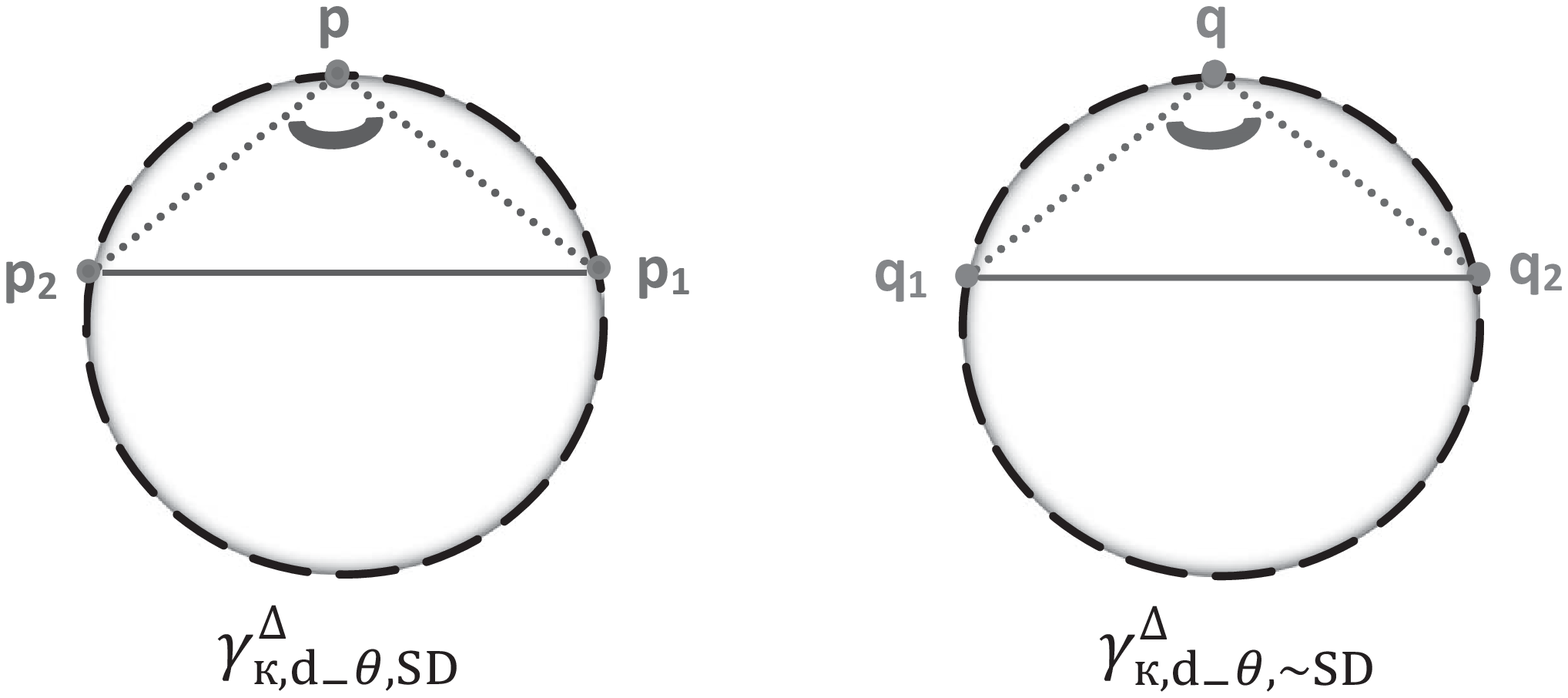} 
\vspace{-6mm}
\caption{4 classes of 3-point equally spaced ordinary meshes where $\mathrm{\mathsf{d\_}}$ is less than the diameter.}
\vspace{-4mm}
\end{figure*} 

\textbf{Corollary 4.13.} \textit{There exists a unique congruent class of 3-point equally spaced ordinary meshes with the same $\mathrm{\kappa}$, $\mathrm{\mathsf{d\_}}$, and angle-type $\mathrm{\theta}$ which are also in the same signature-direction}.

\textbf{Theorem 4.14, Signature-inverse Theorem.} \textit{Equally spaced ordinary meshes $\mathrm{\gamma^{\vartriangle}}$ and $\mathrm{\tilde{\gamma}^{\vartriangle}}$ in $\mathrm{M^{S\Bbb{E}}}$ with the same signature-direction and angle-type which also have identical 2-point signatures in terms of (2) are congruent}.

\textbf{Proof.} According to (2) and just like the Theorem 4.9 we have
\vspace{-1mm}
\begin{eqnarray}
\mathrm{\frac{4\Delta}{abc}(\tilde{p}_{i})=\frac{4\Delta}{abc}(p_{i}); \hspace{5mm} 1<i<n \hspace{5mm} \mbox{and}}
\end{eqnarray}
\vspace{-10mm}
\begin{eqnarray}
\vspace{-7mm}
\mathrm{\tilde{d}_{i-1,i+1}=d_{i-1,i+1}; \hspace{5mm} 2<i<n-1}.
\end{eqnarray}

\vspace{-3mm}
\noindent Since, additionally, $\mathrm{\gamma^{\vartriangle}}$ and $\mathrm{\tilde{\gamma}^{\vartriangle}}$ are equally spaced of the same angle-type and are in the same signature-direction, by Corollary 4.13
\vspace{-3mm}
\begin{eqnarray}
\mathrm{\exists \hspace{.5mm}g_{i}\in S\Bbb{E}(2) \hspace{4mm} \mbox{s.t.} \hspace{4mm} g_{i}: N^{1}_{i}\longmapsto \tilde{N}^{1}_{i}; \hspace{5mm} 2<i<n-1}.
\end{eqnarray} 

\vspace{-3mm}
\noindent To prove that $\mathrm{g_{i}}$s are equal, consider the one-neighborhoods of two successive points $\mathrm{p_{i}}$ and $\mathrm{p_{i+1}}$. From (14) 

\vspace{-4mm}
$\mathrm{\hspace{2.7cm} g_{i}\cdot [p_{i}, p_{i+1}]\longmapsto [\tilde{p}_{i}, \tilde{p}_{i+1}]} \hspace{5mm} \mbox{and} \hspace{5mm} \mathrm{g_{i+1}\cdot [p_{i}, p_{i+1}]\longmapsto [\tilde{p}_{i}, \tilde{p}_{i+1}]}$.
\clearpage

\noindent Hence, according to Lemma 4.3, $\mathrm{g_{i}=g_{i+1}}$ for $\mathrm{2<i<n-1}$. In other words
\vspace{-4mm}
\begin{eqnarray*}
\mathrm{\exists ! \hspace{.7mm} g\in S\Bbb{E}(2) \hspace{4mm} s.t. \hspace{4mm} \tilde{\gamma}^{\vartriangle}=g\cdot \gamma^{\vartriangle}; \hspace{5mm} p_{2}<p_{i}<p_{n-1}}.
\end{eqnarray*}

\vspace{-4mm}
\noindent Finally, Corollary 4.6 extends this result to the whole points of the meshes, meaning $\mathrm{\gamma^{\vartriangle}}$ and $\mathrm{\tilde{\gamma}^{\vartriangle}}$ are $\mathrm{S\Bbb{E}}$-congruent.

We indeed established the following theorem for \textit{fine} meshes as well.

\noindent \hspace{4.4mm}\textbf{Theorem 4.15, Signature-inverse Theorem.} \textit{Two equally spaced ordinary fine meshes $\mathrm{\gamma^{\vartriangle}, \tilde{\gamma}^{\vartriangle}\subset M^{S\Bbb{E}}}$ in the same signature-direction and with identical 2-point JINS, parameterized by (2), are congruent}.
%
%
%==================================================================================================================
%
%

\vspace{2.5mm}
\textbf{$\square \succ$ \textit{Signature-inverse Theorem in terms of (2) - version 2}}
\vspace{2mm}

\textbf{Definition 4.16.} \textit{Let $\mathrm{\gamma^{\vartriangle}=\lbrace p_{i}\rbrace\subset M^{S\Bbb{E}}}$ denote an ordinary mesh. By ``the signed angle" at a point $\mathrm{p_{i}}$ we mean $\mathrm{\vartheta_{i}=SS(p_{i})\cdot \theta_{i}}$ where $\mathrm{\theta_{i}=\hspace{1.2mm}<\hspace{-1.3mm}(p_{i-1}, p_{i}, p_{i+1})}$}. 

\noindent \textit{Also, two ordinary meshes $\mathrm{\gamma^{\vartriangle}}$ and $\mathrm{\tilde{\gamma}^{\vartriangle}}$ have the same ``signed angle-type" if the corresponding angles $\mathrm{\vartheta_{i}}$ and $\mathrm{\tilde{\vartheta}_{i}}$ are both signed-acute, signed-obtuse, or signed-right}.

Definition 4.16, along with Corollary 4.13 and Theorem 4.14, gives the following results.

\textbf{Theorem 4.17.} \textit{There exists a unique congruent class of equally spaced 3-point ordinary meshes with $\mathrm{{\Bbb{E}}}$-curvature $\mathrm{\kappa}$, distance $\mathrm{\mathsf{d\_}}$, and the same signed angle-type $\mathrm{\vartheta}$, see FIG. 12}.
\begin{figure*}[!hbt]
\vspace{-4.5mm}
\includegraphics[angle=0,scale=0.4]{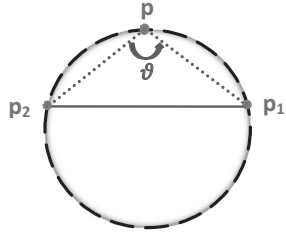} 
\vspace{-6mm}
\caption{The unique class of equally spaced 3-point meshes with the same $\mathrm{\kappa}$, $\mathrm{\mathsf{d\_}}$, and $\mathrm{\vartheta}$.}
\vspace{-4mm}
\end{figure*}  

\textbf{Theorem 4.18, Signature-inverse Theorem.} \textit{Equally spaced ordinary meshes $\mathrm{\gamma^{\vartriangle}}$ and $\mathrm{\tilde{\gamma}^{\vartriangle}}$ in $\mathrm{M^{S\Bbb{E}}}$ with the same signed angle-type and 2-point signature $\mathrm{\Xi_{S\Bbb{E}}^{\vartriangle}(\tilde\gamma^{\vartriangle})=\Xi_{S\Bbb{E}}^{\vartriangle}(\gamma^{\vartriangle})}$ obtained from (2) are congruent}.

Now, we rewrite Theorem 4.8 in terms of the signed angles as follows.

\textbf{Theorem 4.19.} \textit{There is a unique congruent class of 3-point equally spaced meshes with Euclidean curvature $\kappa$ and signed-angle $\mathrm{0<\vartheta_{0}<\pi}$}.

Theorem 4.19 establishes the following result.

\textbf{Theorem 4.20.} \textit{Equally spaced ordinary meshes $\mathrm{\gamma^{\vartriangle}}$ and $\mathrm{\tilde{\gamma}^{\vartriangle}}$ in $\mathrm{M^{S\Bbb{E}}}$ with identical signed angles $\mathrm{0<\tilde\vartheta_{0,i}=\vartheta_{0,i}<\pi}$ and 2-point Euclidean curvatures $\mathrm{\kappa_{\Bbb{E}}^{\vartriangle} ({\tilde\gamma^{\vartriangle}})=\kappa_{\Bbb{E}}^{\vartriangle} ({\gamma^{\vartriangle}})}$ are congruent.}
%
%
%-----------------------------------------------------------------------------------======================================
%
%
\vspace{-5mm}
\subsection{4.2 Unequally spaced meshes in $\mathrm{M^{S\Bbb{E}}}$}
\vspace{-3mm}

\textbf{Definition 4.21 \cite{agh}.} \textit{Let $\mathrm{\gamma^{\vartriangle}=\lbrace p_{i}\rbrace \subset M^{S\Bbb{E}}}$ be an ordinary mesh}. 

\noindent a) The $\mathrm{(m_{1}, m_{2})}$-neighborhood \textit{of a point $\mathrm{p_{i}\in \gamma^{\vartriangle}}$ means $\mathrm{p_{i-m_{1}}, p_{i}, p_{i+m_{2}}\in E}$}.

\noindent b) The 2-point $\mathrm{\Bbb{E}_{m_{1}}^{m_{2}}}$-curvature \textit{of $\mathrm{\gamma^{\vartriangle}}$ is the real-valued function $(m_{1},m_{2})$-$\mathrm{\kappa^{\vartriangle}_{\Bbb{E}}:\gamma^{\vartriangle}\longrightarrow \Bbb{R}}$, given by Definition 2.1 rewritten in the $(m_{1}, m_{2})$-neighborhood of each point}. 

\noindent c) $(m_{1},m_{2})$-$\mathrm{\Xi^{\vartriangle}_{S\Bbb{E}}(\gamma^{\vartriangle})}$ \textit{denotes} the 2-point $\mathrm{S\Bbb{E}_{m_{1}}^{m_{2}}}$-signature \textit{of $\mathrm{\gamma^{\vartriangle}}$ parameterized by (1)-(4) where they are rewritten in the $(m_{1}, m_{2})$-neighborhood of each point in terms of the} $(m_{1},m_{2})$-$\mathrm{\kappa^{\vartriangle}_{\Bbb{E}}}$. 

\noindent d) The $\mathrm{(m_{1},m_{2})}$-angle \textit{at $\mathrm{p_{i}}$ refers to $(m_{1}, m_{2})$-$\mathrm{\theta_{i}=\hspace{1.2mm}<\hspace{-1.3mm}(p_{i-m_{1}}, p_{i}, p_{i+m_{2}})}$ along $\mathrm{\gamma^{\vartriangle}}$. In addition,} the signed $\mathrm{(m_{1},m_{2})}$-angle \textit{at $\mathrm{p_{i}}$ means} $(m_{1},m_{2})$-$\mathrm{\vartheta_{i}=SS(p_{i})\cdot} (m_{1}, m_{2})$-$\mathrm{\theta_{i}}$.

\noindent Moreover, if $\mathrm{m_{1}=m_{2}=m}$ then ``$\mathrm{(m_{1}, m_{2})}$-$\mathrm{\ast}$" will be denoted by ``$\mathrm{\ast^{m}}$".

According to Aghayan \cite{agh1, agh3}, we have the following theorem.

\textbf{Theorem 4.22, Signature Theorem in $\mathrm{M^{S\Bbb{E}}}$ \cite{agh, agh1, agh2}.} \textit{Congruent ordinary meshes have the same (orientation-invariant) 2-point $\mathrm{S\Bbb{E}_{m_{1}}^{m_{2}}}$-signature}.

%
%---------------------------------------------------------------------------------------------------------------
%
\vspace{2mm}
\textbf{$\square \succ$ \textit{Signature-inverse Theorem in terms of (3) and (4)}}
\vspace{1.5mm}

In an unequally spaced 3-point mesh, let $\mathrm{\mathsf{d_{1}}, \mathsf{d_{2}}}$ and $\mathrm{\mathsf{d\_}}$ denote, respectively, the $\mathrm{\Bbb{E}}$-distances between the first and middle points, the end and middle points, and the first and end points. Besides for an unequally spaced n-point mesh, let $\mathrm{(m_{1},m_{2})}$-$\mathrm{\mathsf{d_{1,i}}}$, $\mathrm{(m_{1},m_{2})}$-$\mathrm{\mathsf{d_{2,i}}}$ and $\mathrm{(m_{1},m_{2})}$-$\mathrm{\mathsf{d_{\_i}}}$ denote, respectively, the $\mathrm{\Bbb{E}}$-distances between the first and middle points, the end and middle points, and the first and end points of the $\mathrm{(m_{1},m_{2})}$-neighborhood of the $\mathrm{i^{th}}$ point (for open meshes $\mathrm{ m_{1}<i<n-m_{2}+1}$ where $\mathrm{n>3}$). We also suppose that $\mathrm{(m_{1},m_{2})}$-$\mathrm{\mathsf{d_{2}}>(m_{1},m_{2})}$-$\mathrm{\mathsf{d_{1}}}$.

\textbf{Lemma 4.23.} \textit{There are at most 4 congruent classes of 3-point unequally spaced meshes with $\mathrm{\Bbb{E}}$-curvature $\kappa$ and the side lengths $\mathrm{\mathsf{d}_{1}}$ and $\mathrm{\mathsf{d}_{2}}$ classified as follows.}

\noindent \hspace{1.5mm} a) \textit{If} $\mathrm{\kappa\cdot d_{2}=2}$;$\mathrm{^{1}}$\footnotetext[1]{meaning that the bigger side $\mathrm{\mathsf{d_{2}}}$ is the diameter $\mathrm{2\mathsf{R}}$ of the circumcircle \vspace{-4mm}} \textit{there exist 2 equivalent classes}: \textit{one in $\mathrm{SD}$ denoted by $\mathrm{\gamma^{\vartriangle}_{\kappa, d_{1}, SD}}$ and another in $\mathrm{\sim SD}$ denoted by $\mathrm{\gamma^{\vartriangle}_{\kappa, d_{1}, \sim SD}}$, see FIG. 13.} 
\begin{figure*}[!hbt]
\vspace{-5mm}
\includegraphics[angle=0,scale=0.39]{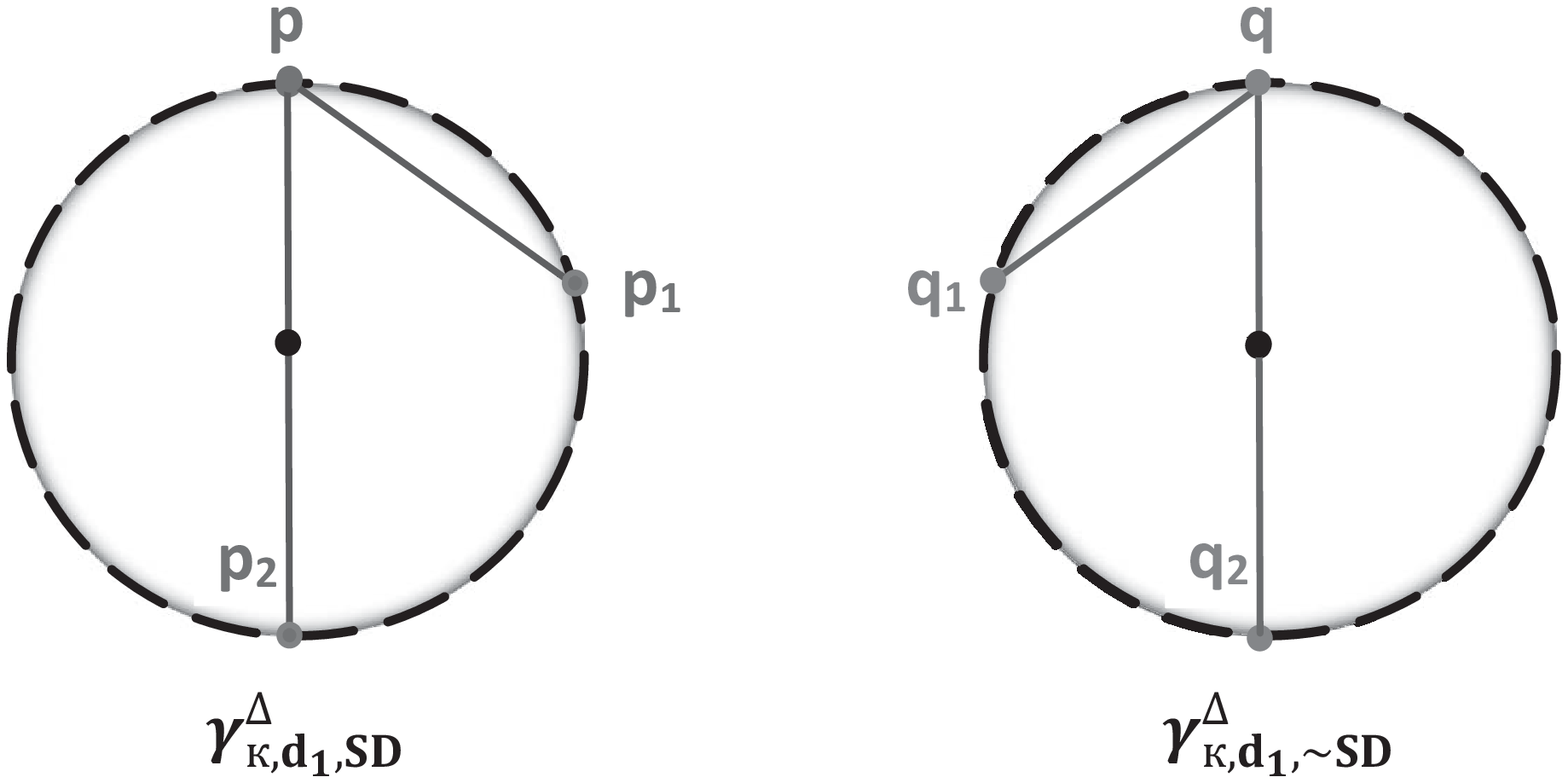}
\vspace{-5mm}
\caption{There is no motion in $\mathrm{S\Bbb{E}}$ to map one to the other.}
\vspace{-3.5cm}
\end{figure*}  
\clearpage

\noindent \hspace{1.5mm} b) \textit{If} $\mathrm{\kappa\cdot \mathsf{d}_{2}\neq 2}$; \textit{there exist 4 congruent classes}: \textit{two in $\mathrm{SD}$ denoted by $\mathrm{\lbrace\gamma^{\vartriangle,1}_{\kappa, \mathsf{d}_{1}, \mathsf{d}_{2}, SD}, \gamma^{\vartriangle,2}_{\kappa, \mathsf{d_{1}}, \mathsf{d_{2}}, SD}\rbrace}$ and the other two in $\mathrm{\sim SD}$ denoted by $\mathrm{\lbrace\gamma^{\vartriangle,3}_{\kappa, \mathsf{d_{1}}, \mathsf{d}_{2}, \sim SD}, \gamma^{\vartriangle,4}_{\kappa, \mathsf{d}_{1}, \mathsf{d}_{2}, \sim SD}\rbrace}$, see FIG. 14.}
\begin{figure*}[!hbt]
\vspace{-4mm}
\includegraphics[angle=0,scale=0.64]{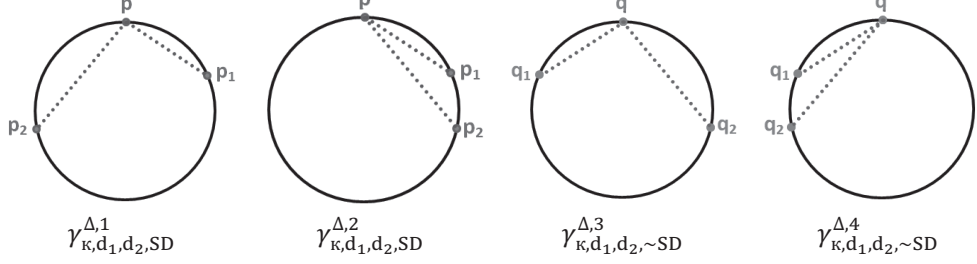}
\vspace{-5.5mm}
\caption{There is no $\mathrm{S\Bbb{E}}$-transformation to coincide one with the other.}
\vspace{-5mm}
\end{figure*} 

We, therefore, have the following theorem.

\textbf{Theorem 4.24.} \textit{In $\mathrm{M^{S\Bbb{E}}}$, there exists a unique congruent class of unequally spaced 3-point ordinary meshes with the following conditions.}

\noindent a) \textit{With the same $\mathrm{\kappa}$, $\mathrm{\mathsf{d}_{1} (\mathsf{d}_{2})}$, and signed angle $\mathrm{\vartheta_{0}}$, denoted by $\mathrm{\gamma^{\vartriangle}_{\kappa, \mathsf{d}_{1}, \vartheta_{0}}}$ $\mathrm{(\gamma^{\vartriangle}_{\kappa, \mathsf{d}_{2}, \vartheta_{0}})}$, see FIG. 15.}
\begin{figure*}[!hbt]
\vspace{-4mm}
\includegraphics[angle=0,scale=0.37]{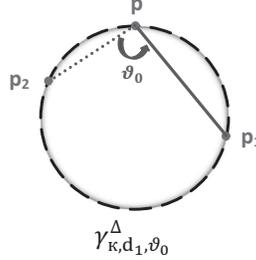}
\vspace{-5mm}
\caption{A unique class with given $\mathrm{\kappa}$, $\mathrm{\mathsf{d}_{1}}$, and signed angle $\mathrm{\vartheta_{0}}$.}
\vspace{-5mm}
\end{figure*} 

\noindent b) \textit{With the same $\mathrm{\kappa}$, $\mathrm{\mathsf{d}_{1} (\mathsf{d}_{2})}$, and signed angle-type $\vartheta \geq 90^{o}$, denoted by $\mathrm{\gamma^{\vartriangle}_{\kappa, \mathsf{d}_{1}, \vartheta\geq 90^{o}}}$ $\mathrm{(\gamma^{\vartriangle}_{\kappa, \mathsf{d}_{2}, \vartheta\geq 90^{o}})}$}.

\noindent c) \textit{With the same $\mathrm{\kappa}$, $\mathrm{\mathsf{d}_{1} (\mathsf{d}_{2})}$, $\mathrm{\mathsf{d}\_}$, and signed angle-type $\vartheta$, denoted by $\mathrm{\gamma^{\vartriangle}_{\kappa, \mathsf{d}_{1}, \mathsf{d}\_, \vartheta}}$ $\mathrm{(\gamma^{\vartriangle}_{\kappa, \mathsf{d}_{2}, \mathsf{d}\_, \vartheta})}$, see FIG. 16}.
\begin{figure*}[!hbt]
\vspace{-8mm}
\includegraphics[angle=0,scale=0.375]{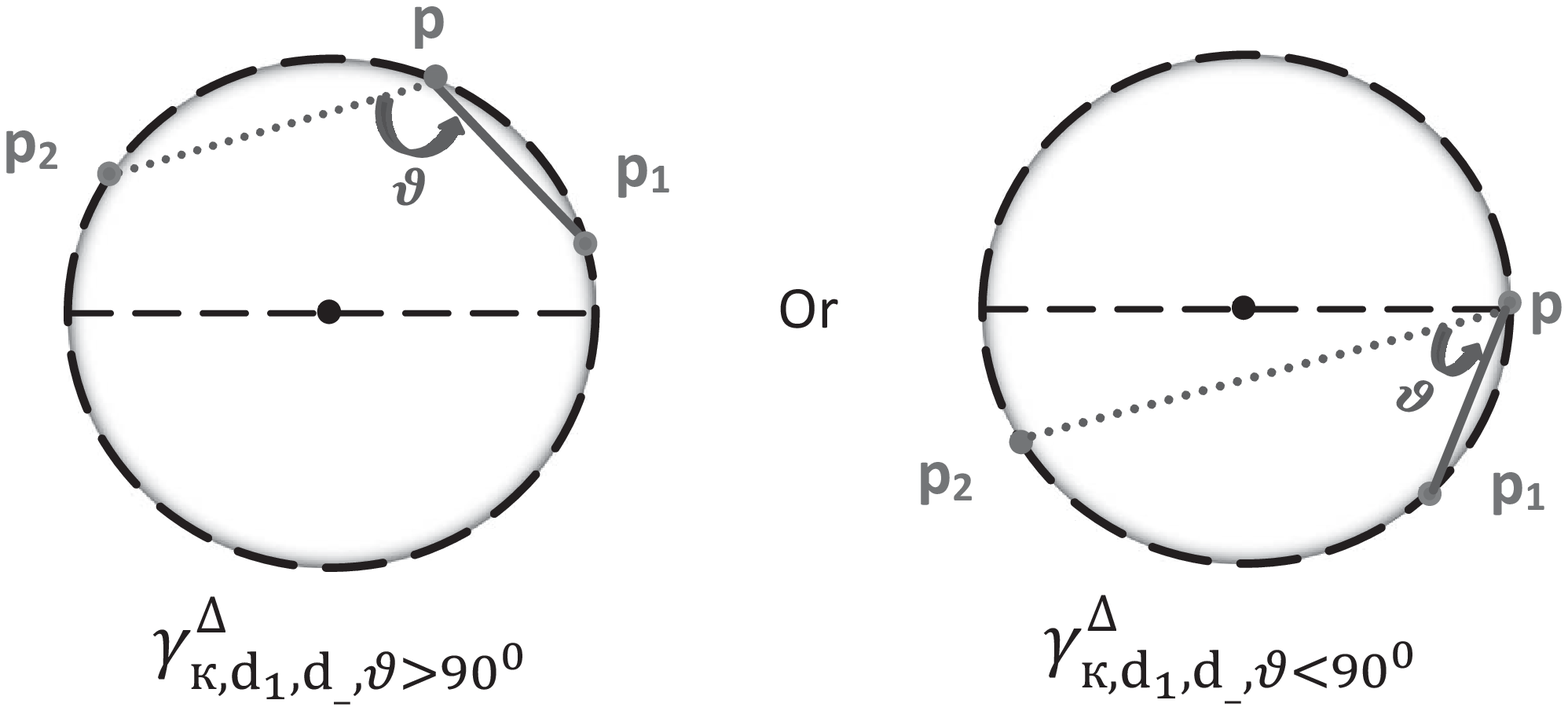} 
\vspace{-5mm}
\caption{A unique class with given $\mathrm{\kappa}$, $\mathrm{\mathsf{d}_{1}}$, $\mathrm{\mathsf{d}\_}$, and signed angle-type.}
\vspace{-5mm}
\end{figure*} 
%
%-----------------------------------------------------------------------------------------------------------------------------------------
%

\textbf{Theorem 4.25, Signature-inverse Theorem.} \textit{Unequally spaced ordinary meshes $\mathrm{\gamma^{\vartriangle}}$ and $\mathrm{\tilde{\gamma}^{\vartriangle}}$ in $\mathrm{M^{S\Bbb{E}}}$ with identical lengths $d_{i-1,i+1}$, the same signed (1,2)-angle type, and the same 2-point $\mathrm{S\Bbb{E}_{1}^{2}}$-signature in terms of (3) are congruent}.

\textbf{Proof.} According to the definition of $\mathrm{S\Bbb{E}_{1}^{2}}$-signatures in terms of (3)
\[\left\{\begin{array}{cl}
\mathrm{(1,2)\mbox{-}\kappa_{\Bbb{E}}^{\vartriangle}(\tilde{p}_{i})=(1,2)\mbox{-}\kappa_{\Bbb{E}}^{\vartriangle}(p_{i})}\hspace{6.3cm} \mbox{and}\\ [1.5mm]
\mathrm{\displaystyle \frac{(1,2)\mbox{-}\kappa_{\Bbb{E}}^{\vartriangle}(\tilde{p}_{i+1})-(1,2)\mbox{-}\kappa_{\Bbb{E}}^{\vartriangle}(\tilde{p}_{i-1})}{\tilde{d}_{i-1, i+2}}=\frac{(1,2)\mbox{-}\kappa_{\Bbb{E}}^{\vartriangle}(p_{i+1})-(1,2)\mbox{-}\kappa_{\Bbb{E}}^{\vartriangle}(p_{i-1})}{d_{i-1, i+2}}}\end{array}\right.\]
\clearpage

\noindent where $\mathrm{1<i<n-2}$. Therefore
\vspace{-3mm}
\begin{eqnarray}
\mathrm{(1,2)\mbox{-}\kappa_{\Bbb{E}}^{\vartriangle}(\tilde{p}_{i})=(1,2)\mbox{-}\kappa_{\Bbb{E}}^{\vartriangle}(p_{i})=\kappa_{i}; \hspace{5mm} 1<i<n-1, \hspace{5mm} \mbox{and}}
\end{eqnarray}
\vspace{-13mm}
\begin{eqnarray}
\mathrm{\tilde{d}_{i-1,i+2}=d_{i-1,i+2}=(1,2)\mbox{-}\mathsf{d }\_{i}; \hspace{5mm} 1<i<n-2}.
\end{eqnarray}

\vspace{-3mm}
\noindent Additionally, $\mathrm{(1,2)\mbox{-}d_{2,i}=(1,2)\mbox{-}\tilde{d}_{2,i}}$ and the signed angles $\mathrm{(1,2)\mbox{-}\vartheta_{i}}$ and $\mathrm{(1,2)\mbox{-}\tilde{\vartheta}_{i}}$ have the same type, therefore, according to Theorem 4.24-part c)
\vspace{-3mm}
\begin{eqnarray}
\mathrm{\exists \hspace{.5mm}g_{i}\in S\Bbb{E}(2) \hspace{4mm} \mbox{s.t.} \hspace{4mm} g_{i}: (1,2)\mbox{-}N_{i}\longmapsto (1,2)\mbox{-}\tilde{N}_{i}; \hspace{5mm} 1<i<n-2}.
\end{eqnarray}

\vspace{-3mm}
\noindent Now, consider the corresponding triangles $\mathrm{\bigtriangleup(p_{i-1}, p_{i}, p_{i+1})}$ and $\mathrm{\bigtriangleup(\tilde{p}_{i-1}, \tilde{p}_{i}, \tilde{p}_{i+1})}$. From (17)

$\mathrm{\hspace{4.35cm} d_{1,i}=\tilde{d}_{1,i}, \hspace{4mm} d_{2,i}=\tilde{d}_{2,i}, \hspace{4mm} \mbox{and} \hspace{5mm} d\__{i}=\tilde{d}\__{i}}$

\noindent which, with regard to the Lemma 4.3, results in $\mathrm{g_{i}=g_{i+1}}$ for $\mathrm{1<i<n-2}$. In other words

$\mathrm{\hspace{3.5cm} \exists ! \hspace{.7mm} g\in S\Bbb{E}(2) \hspace{4mm} s.t. \hspace{4mm} \tilde{\gamma}^{\vartriangle}=g\cdot \gamma^{\vartriangle}; \hspace{5mm} 1<i<n-2}$

\noindent showing the congruence of $\mathrm{\gamma^{\vartriangle}}$ and $\mathrm{\tilde{\gamma}^{\vartriangle}}$ where they are closed. It also proves that for the open meshes the congruence satisfies where $\mathrm{p_{1}\leq p_{i}<p_{n}}$. In this case, if $\mathrm{\gamma^{\vartriangle}}$ and $\mathrm{\tilde{\gamma}^{\vartriangle}}$ have the same $\mathrm{(1,2)}$-$\mathrm{d\__{(n-2)}}$, the result extends to their whole points.
%
%Moreover, by considering the corresponding triangles $\mathrm{\bigtriangleup(p_{n-4}, p_{n-3}, p_{n-2})}$ and $\mathrm{\bigtriangleup(\tilde{p}_{n-4}, \tilde{p}_{n-3}, \tilde{p}_{n-2})}$ we can extend (26) to $\mathrm{i<n-1}$. And therefore, $\mathrm{\gamma^{\vartriangle}}$ and $\mathrm{\tilde{\gamma}^{\Delta}}$ are $\mathrm{S\Bbb{E}}$-congruent for all points of the given curves if they are closed, and for $\mathrm{p_{1}\leq p_{i}<p_{n-1}}$ if they are open.

%
%-----------------------------------------------------------------------------------------------------------------------------------------
%
%

\textbf{Theorem 4.26, Signature-inverse Theorem.} \textit{Let $\mathrm{\gamma^{\vartriangle}, \tilde{\gamma}^{\vartriangle}\subset M^{S\Bbb{E}}}$ denote two unequally spaced ordinary meshes with the same $\mathrm{\Bbb{E}^{1}_{3}}$-curvature and signed angle $\mathrm{(3,1)}$-$\mathrm{\vartheta_{0,i}}$. Let also they have the same signed 3-angle type, identical lengths $d^{3}_{1,i}$, and identical 2-point $\mathrm{S\Bbb{E}^{3}}$-signatures in terms of (4). Then, $\mathrm{\gamma^{\vartriangle}}$ and $\mathrm{\tilde{\gamma}^{\vartriangle}}$ are congruent}.

\textbf{Proof.} With regard to the definition of $\mathrm{S\Bbb{E}^3}$-signatures in terms of (4)
\vspace{-2.5mm}
\[\left\{\begin{array}{cl}
\mathrm{\kappa_{\Bbb{E}}^{\vartriangle,3}(\tilde{p}_{i})=\kappa_{\Bbb{E}}^{\vartriangle,3}(p_{i})}\hspace{4.7cm} \mbox{and}\\ [1.5mm]
\mathrm{\displaystyle \frac{\kappa_{\Bbb{E}}^{\vartriangle,3}(\tilde{p}_{i+1})-\kappa_{\Bbb{E}}^{\vartriangle,3}(\tilde{p}_{i-1})}{\tilde{d}_{i-3, i+3}}=\frac{\kappa_{\Bbb{E}}^{\vartriangle,3}(p_{i+1})-\kappa_{\Bbb{E}}^{\vartriangle,3}(p_{i-1})}{d_{i-3, i+3}}}\end{array}\right.\]

\vspace{-2.5mm}
\noindent where $\mathrm{4<i<n-3}$. As a result
\vspace{-3.5mm}
\begin{eqnarray}
\mathrm{\kappa_{\Bbb{E}}^{\vartriangle,3}(\tilde{p}_{i})=\kappa_{\Bbb{E}}^{\vartriangle,3}(p_{i})=\kappa_{i}; \hspace{5mm} 3<i<n-2 \hspace{5mm} \mbox{and}}
\end{eqnarray}
\vspace{-13mm}
\begin{eqnarray}
\mathrm{\tilde{d}_{i-3,i+3}=d_{i-3,i+3}=d^{3}_{\_i}; \hspace{5mm} 4<i<n-3}.
\end{eqnarray}

\vspace{-3mm}
\noindent Besides, $\mathrm{d^{3}_{1,i}=\tilde{d}^{3}_{1,i}}$ and the signed angles $\mathrm{\tilde{\vartheta}_{i}^{3}}$ and $\mathrm{\vartheta_{i}^{3}}$ have the same type, therefore, according to Theorem 4.24-part c)
\vspace{-4.5mm}
\begin{eqnarray}
\mathrm{\exists g_{i}\in S\Bbb{E}(2) \hspace{4mm} \mbox{s.t.} \hspace{4mm} g_{i}: N^{3}_{i}\longmapsto \tilde{N}^{3}_{i}; \hspace{5mm} 4<i<n-3}.
\end{eqnarray} 

\vspace{-3.5mm}
Now, through the following steps, we prove that $\mathrm{g_{i}}$s are identical.

\noindent Step 1. 

By (20)
\vspace{-7mm}
\begin{eqnarray*}
\mathrm{g_{i}\cdot [p_{i}, p_{i+3}]\longmapsto [\tilde{p}_{i}, \tilde{p}_{i+3}]} \hspace{5mm} \mbox{and} \hspace{5mm} \mathrm{g_{i+3}\cdot [p_{i}, p_{i+3}]\longmapsto [\tilde{p}_{i}, \tilde{p}_{i+3}]}
\end{eqnarray*} 

\vspace{-3mm}
\noindent which, according to Lemma 4.3, result in $\mathrm{g_{i}=g_{i+3}}$ and thus there would be only the following three group transformations:
\clearpage
\begin{equation}
\begin{cases}
\vspace{-1mm}
\mathrm{g_{0}:N^{3}_{i}\longmapsto \tilde{N}^{3}_{i}} \hspace{6mm} & \text{where } \hspace{2.5mm} \mathrm{i=0 \hspace{1mm} mod 3},
\vspace{-1mm}
\\
\vspace{-1mm}
\mathrm{g_{1}:N^{3}_{i}\longmapsto \tilde{N}^{3}_{i}} \hspace{6mm} & \text{where } \hspace{2.5mm} \mathrm{i=1 \hspace{1mm} mod 3},
\vspace{-1mm}
\\
\vspace{-1mm}
\mathrm{g_{2}:N^{3}_{i}\longmapsto \tilde{N}^{3}_{i}} \hspace{6mm} & \text{where } \hspace{2.5mm} \mathrm{i=2 \hspace{1mm} mod 3}.
\end{cases}
\end{equation}

\vspace{-1mm}
\noindent which we need to show they are equal.

\noindent Step 2. 

Let $\mathrm{i=0 \hspace{2mm} mod 3}$. From (20), we have $\mathrm{(3,1)\mbox{-}d_{1,i}=(3,1)\mbox{-}\tilde{d}_{1,i}}$ and additionally

$\mathrm{\hspace{2.1cm} (3,1)}$-$\mathrm{\kappa^{\vartriangle}_{\Bbb{E}}(p_{i})=}$ $\mathrm{(3,1)}$-$\mathrm{\kappa^{\vartriangle}_{\Bbb{E}}(\tilde{p}_{i}) \hspace{5mm} \mbox{and} \hspace{6mm} (3,1)}$-$\mathrm{\vartheta_{0,i}=}$ $\mathrm{(3,1)}$-$\mathrm{\tilde{\vartheta}_{0,i}}$. 

\noindent Therefore, with regard to Theorem 4.24-part a)
\vspace{-2.5mm}
\begin{eqnarray}
\mathrm{\exists g^{0}_{i}\in S\Bbb{E}(2) \hspace{4mm} \mbox{s.t.} \hspace{4mm} g^{0}_{i}: (3,1)\mbox{-}N_{i}\longmapsto (3,1)\mbox{-}\tilde{N}_{i}; \hspace{5mm} 4<i<n-2}.
\end{eqnarray} 

\vspace{-2.5mm}
\noindent Identity (22), along with (20), indicates that
\vspace{-3mm}
\begin{eqnarray*}
\mathrm{g_{0}\cdot [p_{i-3}, p_{i}]\longmapsto [\tilde{p}_{i-3}, \tilde{p}_{i}]} \hspace{5mm} \mbox{and} \hspace{5mm} \mathrm{g^{0}_{i}\cdot [p_{i-3}, p_{i}]\longmapsto [\tilde{p}_{i-3}, \tilde{p}_{i}]}
\end{eqnarray*} 

\vspace{-2.5mm}
\noindent which, according to Lemma 4.3, results in $\mathrm{g^{0}_{i}=g_{0}}$. Hence, by (22)
\vspace{-1mm}
\begin{eqnarray}
\mathrm{\exists \hspace{.5mm}g_{0}: (3,1)\mbox{-}N_{i}\longmapsto (3,1)\mbox{-}\tilde{N}_{i}; \hspace{5mm} 4<i<n-2}.
\end{eqnarray} 

\vspace{-5mm}
\noindent Step 3.

Consider now the triangles $\mathrm{\bigtriangleup(p_{i}, p_{i+1}, p_{i+3})}$ and $\mathrm{\bigtriangleup(\tilde{p}_{i}, \tilde{p}_{i+1}, \tilde{p}_{i+3})}$. From (20) and (23)
\vspace{-3mm}
\begin{eqnarray*}
\mathrm{d_{i, i+3}=\tilde{d}_{i, i+3}}, \hspace{3mm} \mathrm{d_{i}=\tilde{d}_{i}}, \hspace{3mm} \mbox{and} \hspace{4mm} \mathrm{<\hspace{-1.3mm}(p_{i+1}, p_{i}, p_{i+3}) =<\hspace{-1.3mm}(\tilde{p}_{i+1}, \tilde{p}_{i}, \tilde{p}_{i+3})}.
\end{eqnarray*} 

\vspace{-3mm}
\noindent Accordingly, (1,2)-triangles at $\mathrm{p_{i+1}\in \gamma^{\vartriangle}}$ and $\mathrm{\tilde{p}_{i+1}\in \tilde{\gamma}^{\vartriangle}}$ are congruent, where $\mathrm{4<i<n-2}$, resulting
\vspace{-6mm}
\begin{eqnarray}
\mathrm{g_{0}: (1,2)\mbox{-}N_{i+1}\longmapsto (1,2)\mbox{-}\tilde{N}_{i+1}}
\end{eqnarray} 

\vspace{-6mm}
\noindent and
\vspace{-6mm}
\begin{eqnarray}
\mathrm{g_{0}: (1,2)\mbox{-}N_{i-2}\longmapsto (1,2)\mbox{-}\tilde{N}_{i-2}}
\end{eqnarray} 

\vspace{-3mm}
\noindent which, along with (23), result in

$\mathrm{\hspace{1.2cm} (2,1)\mbox{-}d_{1,i}=(2,1)\mbox{-}\tilde{d}_{1,i}, \hspace{3mm} (2,1)\mbox{-}\vartheta_{i}=(2,1)\mbox{-}\tilde{\vartheta}_{i}, \hspace{3mm} \mbox{and} \hspace{4mm} (2,1)\mbox{-}d_{2,i}=(2,1)\mbox{-}\tilde{d}_{2,i}}$.

\noindent In other words, we have
\vspace{-3mm}
\begin{eqnarray}
\mathrm{g_{0}: (2,1)\mbox{-}N_{i}\longmapsto (2,1)\mbox{-}\tilde{N}_{i}; \hspace{5mm} 6<i<n-3}.
\end{eqnarray}

\vspace{-4mm}
Step 4. 

From (23) and (26)

$\mathrm{\hspace{1.8cm} g_{1}\cdot [p_{i-2}, p_{i+1}]\longmapsto [\tilde{p}_{i-2}, \tilde{p}_{i+1}] \hspace{4mm} \mbox{and} \hspace{4.5mm} g_{0}\cdot [p_{i-2}, p_{i+1}]\longmapsto [\tilde{p}_{i-2}, \tilde{p}_{i+1}]}$

\noindent which, according to Lemma 4.3, result in $\mathrm{g_{1}=g_{0}}$. 

Similarly, we can prove $\mathrm{g_{2}=g_{0}}$. Thus

$\mathrm{\hspace{2.6cm} \exists ! \hspace{.7mm} g\in S\Bbb{E}(2) \hspace{4mm} s.t. \hspace{4mm} \tilde{\gamma}^{\vartriangle}=g\cdot \gamma^{\vartriangle} \hspace{5mm} \mbox{for} \hspace{3mm} 6<i<n-3}$

\noindent showing the congruence of $\mathrm{\gamma^{\vartriangle}}$ and $\mathrm{\tilde{\gamma}^{\vartriangle}}$ where they are closed. It also proves that for the open meshes the congruence satisfies where $\mathrm{p_{2}\leq p_{i}<p_{n-1}}$. In this case, adding extra conditions for the first and end points in the open case, for example $\mathrm{\vartheta^{3}_{4} \geq 90^{o}}$ or $\mathrm{\gamma^{\vartriangle}}$ and $\mathrm{\tilde{\gamma}^{\vartriangle}}$ have identical corresponding 3-angles $\mathrm{\vartheta^{3}_{0,4}}$ and $\mathrm{\vartheta^{3}_{0,n-3}}$, extends the result to the whole points.
\clearpage
%
%-----------------------------------------------------------------------------------------------------------------------------------------
%
\vspace{-5mm}
\subsection{4.3 The Host Theorem and Signature-inverse Theorem}
\vspace{-4mm}
%
%-----------------------------------------------------------------------------------------------------------------------------------------
%

40 guests are supposed to sit around a dinner table. The host has asked the party planner to provide 3, 4, or 5 different types of rare wildflowers and put them on the table in front of each guest between the dinner and dessert so that the following conditions must be met. 

\noindent $\surd$ Depending on the variety in the selected flowers (2, 3, or 5), every second, third, or fourth person would receive the same type of flower. 

\noindent $\surd$ The party planner must not look at the guests. He should just step forward and put each flower on the table with his right hand. 

\noindent $\surd$ In addition, each gust must have one and only one flower.

Below we bring a general mathematical setting for this problem.

\textbf{The Host Problem.} \textit{Let $\mathrm{\lbrace p_{i}\rbrace^{n}_{1}\subset \Bbb{R}^2}$ be a closed ordinary mesh and 1$<$m$<$n. We start from $\mathrm{p_{1}}$ and connect every $m^{th}$ point by a straight line until getting back to $\mathrm{p_{1}}$. The question is for what ``m" we meet all points of the mesh uniquely on our way.}

According to the group of primitive residue classes modulo \textit{m} in modular arithmetic, we bring in the following theorem. 

\textbf{The Host Theorem.} \textit{In The Host Problem, the whole points of the given mesh are met iff after n steps we return to the start point iff $gcd^{1}(m,n)=1$}. \footnotetext[1]{gcd refers to the greatest common divisor.}

We, therefore, have the following corollary.

\textbf{Corollary.} \textit{The party planner can fulfill the task by $\varphi(n)$ ways where $\varphi(n)$ denotes Euler's totient (phi) function.}

\noindent As an example, consider a 10-point mesh $\mathrm{\lbrace p_{i}\rbrace}$ and let $\mathrm{m_{1}=3}$ and $\mathrm{m_{2}=4}$. Since gcd(3,10)=1, every point will be met, see FIG. 17-left, and since gcd$\mathrm{(4,10)\neq 1}$, there is no way one meets all points on their way, see FIG. 17-right.
\begin{figure*}[!hbt]
\vspace{-3.5mm}
\includegraphics[angle=0,scale=0.39]{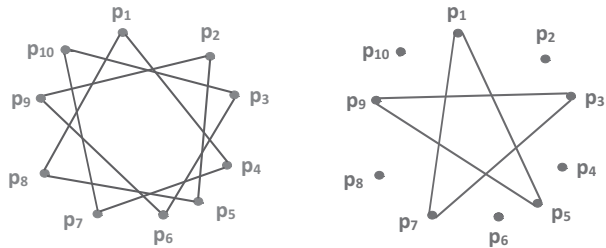}
\vspace{-5.5mm}
\caption{Illustrates The Host Theorem.}
\vspace{-4.5mm}
\end{figure*}  

As a result, since gcd(3,40)=1, the party planner must provide 3 different types of flowers and place them in different baskets behind the host and do his task with no worries.
\clearpage

The Host Theorem results in the following Signature-inverse Theorem for closed meshes.

\noindent \textbf{The Host Theorem and Signature-inverse Theorem.}
\textit{Let $\mathrm{\gamma^{\vartriangle}, \tilde{\gamma}^{\vartriangle}\subset M^{S\Bbb{E}}}$ be two n-point unequally spaced ordinary closed meshes, where $n\neq0 \space\ mod \hspace{.8mm} 3$. Let they have the same signed angle-type $\vartheta_{i}^{3}$, corresponding length $d^{3}_{1,i}$, and 2-point $\mathrm{S\Bbb{E}^{3}}$-signature resulted  from identity (4). Then, $\gamma^{\vartriangle}$ and $\tilde{\gamma}^{\vartriangle}$ are congruent}.
%
%
%================================================================================================================
%
%
\vspace{-4mm}
\section{5. Signature-inverse Theorem in $\mathrm{M^{S\Bbb{A}}}$}
\vspace{-2mm}

\textbf{Definition 5.1.} \emph{Let $\mathrm{\gamma^{\vartriangle}}$ be a convex ordinary mesh. A $\mathrm{p_{i}\in \gamma^{\vartriangle}}$ is in} $\mathrm{\Bbb{A}}$-signature direction $\mathrm{(SD_{\Bbb{A}})}$ \textit{if its two-neighborhood is counterclockwise on the unique conic section passing through it. Otherwise, we say $\mathrm{p_{i}}$ is in $\mathrm{\sim SD_{\Bbb{A}}}$}. 

\noindent \textit{In addition, two n-point convex meshes $\mathrm{\gamma^{\vartriangle}=\lbrace p_{i} \rbrace}$ and $\mathrm{\tilde\gamma^{\vartriangle}=\lbrace \tilde{p}_{i} \rbrace}$ are in the same $\mathrm{\Bbb{A}}$-signature direction iff for $1<i<n$ the corresponding points $\mathrm{p_{i}}$ and $\mathrm{\tilde{p}_{i}}$ are in $\mathrm{SD_{\Bbb{A}}}$ or in $\mathrm{\sim SD_{\Bbb{A}}}$}.

\noindent \hspace{5mm}\textbf{Definition 5.2.} \emph{Let $\mathrm{\gamma^{\vartriangle}=\lbrace p_{i} \rbrace}$ denote a n-point ordinary convex mesh. For $\mathrm{2<i<n-1}$,} the $\mathrm{\Bbb{A}}$-arc length set \textit{of $\mathrm{\gamma^{\vartriangle}}$ at a point $\mathrm{p_{i}}$ is the ordered set $\mathrm{L^{i}=\lbrace L_{k} \mid i-2\leq k\leq i+1\rbrace}$ given by Theorem 2.3 in the two-neighborhood of $\mathrm{p_{i}}$.}

\noindent \textit{Also, ordinary convex meshes $\mathrm{\gamma^{\vartriangle}}$ and $\mathrm{\tilde\gamma^{\vartriangle}}$ have the same $\mathrm{\Bbb{A}}$-arc length set if their corresponding points have the same $\mathrm{\Bbb{A}}$-arc length set.}

\textbf{Definition 5.3.} \textit{Let $\mathrm{\gamma^{\vartriangle}=\lbrace p_{i} \rbrace\subset M^{S\Bbb{A}}}$ be an ordinary convex mesh point}. 

\noindent a) \textit{A positive $\mathrm{\Bbb{A}}$-curvature point $\mathrm{p_{i}}$ has} the fine-area \textit{if the area if its approximating ellipse is not less than the area of the elliptical sector surrounded by the two-neighborhood of $\mathrm{p_{i}}$}.  

\noindent b) \emph{A negative $\mathrm{\Bbb{A}}$-curvature point 
$\mathrm{p_{i}}$} is in the fine-position \textit{if the points in its two-neighborhood are located on the same branch of the approximating hyperbola. In other words, the Euclidean distance between every two successive points of this neighborhood is less than $\mathrm{\mu=2 \sqrt{\frac{-F}{\lambda_{1}^2S}} (p_{i})}$ (twice the semi-major axis of the approximating hyperbola), where F and S are the first and second $\mathrm{\Bbb{A}}$-invariants of $\mathrm{\gamma^{\vartriangle}}$ and $\mathrm{\lambda}$ is the large root of the equation $\mathrm{(\lambda^2-(A+C)\lambda+S)(p_{i})=0}$ in which A and C are the affine-functions of $\gamma^{\vartriangle}$}.

\textbf{Definition 5.4.} \emph{An ordinary convex mesh $\mathrm{\gamma^{\vartriangle}}$ is called} $\mathrm{\Bbb{A}}$-fine \textit{if each point $\mathrm{p_{i}\in \gamma^{\vartriangle}}$ with a positive $\mathrm{\Bbb{A}}$-curvature has the fine-area and with a negative $\mathrm{\Bbb{A}}$-curvature, is in the fine-position.}

It is not hard to demonstrate the following classifications in the affine case.

\textbf{Lemma 5.5} a) \textit{In $\mathrm{M^{S\Bbb{A}}}$, there are 2 congruent classes of five-point ordinary convex meshes with $\mathrm{\Bbb{A}}$-curvature $\mathrm{\kappa_{\Bbb{A}}>0}$, $\mathrm{\Bbb{A}}$-arc length set L, and also their middle points have the fine-area. One in $\mathrm{SD_{\Bbb{A}}}$ denoted by $\mathrm{\gamma^{\vartriangle}_{\kappa_{\Bbb{A}}>0, L, +, SD_{\Bbb{A}}}}$ and another in $\mathrm{\sim SD_{\Bbb{A}}}$ denoted by $\mathrm{\gamma^{\vartriangle}_{\kappa_{\Bbb{A}}>0, L, +, \sim SD_{\Bbb{A}}}}$}.
\clearpage

\noindent b) \textit{In $\mathrm{M^{S\Bbb{A}}}$, there exist 2 congruent classes of five-point ordinary convex meshes with $\mathrm{\kappa_{\Bbb{A}}<0}$, L, and their middle points are in the fine-position. One in $\mathrm{SD_{\Bbb{A}}}$ denoted by $\mathrm{\gamma^{\vartriangle}_{\kappa_{\Bbb{A}}<0, L,+, SD_{\Bbb{A}}}}$ and another in $\mathrm{\sim SD_{\Bbb{A}}}$ denoted by $\mathrm{\gamma^{\vartriangle}_{\kappa_{\Bbb{A}}<0, L, +, \sim SD_{\Bbb{A}}}}$}. 

\noindent c) \textit{In $\mathrm{M^{S\Bbb{A}}}$, there exist 2 congruent classes of five-point ordinary convex meshes with $\mathrm{\kappa_{\Bbb{A}}=0}$ and $\mathrm{\Bbb{A}}$-arc length set L. One in $\mathrm{SD_{\Bbb{A}}}$ and another in $\mathrm{\sim SD_{\Bbb{A}}}$ denoted respectively by $\mathrm{\gamma^{\vartriangle}_{\kappa_{\Bbb{A}}=0, L, SD_{\Bbb{A}}}}$ and $\mathrm{\gamma^{\vartriangle}_{\kappa_{\Bbb{A}}=0, L, \sim SD_{\Bbb{A}}}}$}. 

\textbf{Lemma 5.6 \cite{fel2}.} \textit{The 3-point action $\mathrm{S\Bbb{A}(2)\times M^{3}\longrightarrow M^{3}}$ defined by} 

$\mathrm{\hspace{4cm} g\cdot (p_{1}, p_{2}, p_{3})\longmapsto (g\cdot p_{1}, g\cdot p_{2}, g\cdot p_{3})}$

\noindent \textit{is free on $\mathrm{M^{3}\setminus \mathsf{C}}$, where} $\mathrm{\mathsf{C}=\lbrace (p_{1}, p_{2}, p_{3})\mid \ p_{1}, p_{2},}$ and $\mathrm{p_{3}}$ are collinear $\rbrace$.

\noindent \hspace{3.7mm}Now, we present Signature-inverse Theorem in terms of identity (6).

\noindent \hspace{3.7mm}\textbf{Theorem 5.7, Signature-inverse Theorem.} \textit{Ordinary $\mathrm{\Bbb{A}}$-fine convex meshes $\mathrm{\gamma^{\vartriangle}}$ and $\mathrm{\tilde{\gamma}^{\vartriangle}}$ in $\mathrm{M^{S\Bbb{A}}}$ with never-zero $\Bbb{A}$-curvatures and the same $\mathrm{\Bbb{A}}$-arc length set and 3-point $\mathrm{S\Bbb{A}}$-signature are congruent}.

\textbf{Proof.} According to (6) and just like the Theorem 4.9 
\vspace{-3mm}
\[\left\{\begin{array}{cl}
\mathrm{\frac{S}{F^{2/3}}(\tilde{p}_{i})=\frac{S}{F^{2/3}}(p_{i})}\hspace{4.4cm} \mbox{and}\\ [1.5mm]
\mathrm{\displaystyle \frac{\frac{S}{F^{2/3}}(\tilde{p}_{i+1})-\frac{S}{F^{2/3}}(\tilde{p}_{i-1})}{\tilde{L}_{i-1,i+1}}=\frac{\frac{S}{F^{2/3}}(p_{i+1})-\frac{S}{F^{2/3}}(p_{i-1})}{L_{i-1,i+1}}}\end{array}\right.\] 

\vspace{-3mm}
\noindent where $\mathrm{3<i<n-2}$. Hence
\vspace{-2mm}
\begin{eqnarray}
\mathrm{\frac{S}{F^{2/3}}(\tilde{p}_{i})=\frac{S}{F^{2/3}}(p_{i})=\kappa_{\Bbb{A},i} \hspace{5mm} 2<i<n-1 \hspace{5mm} \mbox{and}}
\end{eqnarray}
\vspace{-12mm}
\begin{eqnarray}
\mathrm{\tilde{L}_{i-1,i+1}=L_{i-1,i+1}; \hspace{5mm} 3<i<n-2}.
\end{eqnarray}

\vspace{-2mm}
\noindent In addition, in the two-neighborhoods of each corresponding points $\mathrm{p_{i}\in \gamma^{\vartriangle}}$ and $\mathrm{\tilde{p}_{i}\in \tilde{\gamma}^{\vartriangle}}$, we have 

\vspace{-4mm}
$\mathrm{\hspace{3.8cm} \tilde{L}_{i-1}=L_{i-1} \hspace{5mm} \mbox{and} \hspace{5mm} \tilde{L}_{i}=L_{i}; \hspace{5mm} 2<i<n-1},$

\noindent which, along with (27), gives

$\mathrm{\hspace{2.7cm} \mbox{Area}(\tilde{p}_{i-1}, \tilde{p}_{i},\tilde{p}_{i+1})=\mbox{Area}(p_{i-1},p_{i},p_{i+1}); \hspace{5mm} 2<i<n-1}.$

\vspace{-2mm}
\noindent Therefore
\vspace{-4mm}
\begin{eqnarray}
\mathrm{\hspace{5mm} \exists \hspace{.5mm} g_{i}\in S\Bbb{A}(2) \hspace{3mm} \mbox{s.t.} \hspace{3mm} g_{i}: \bigtriangleup(\tilde{p}_{i-1},\tilde{p}_{i},\tilde{p}_{i+1})\longrightarrow \bigtriangleup(p_{i-1},p_{i},p_{i+1}); \hspace{5mm} 3<i<n-2}.
\end{eqnarray}

\vspace{-3mm}
\noindent On the other hand, with regard to the hypothesis, Lemma 5.5 indicates that 
\vspace{-3mm}
\begin{eqnarray}
\mathrm{\exists \hspace{.7mm}g^{0}_{i}: N^{2}_{i}\longrightarrow \tilde{N}^{2}_{i}; \hspace{5mm} 2<i<n-1}
\end{eqnarray} 

\vspace{-5mm}
\noindent at any point $\mathrm{p_{i}}$. 

\noindent Identities (29) and (30) along with Lemma 5.6 result in $\mathrm{g_{i-1}=g_{i}=g_{i+1}=g^{0}_{i}}$. Hence
\vspace{-3mm}
\begin{eqnarray*}
\mathrm{\exists ! \hspace{.7mm} g\in S\Bbb{A}(2) \hspace{4mm} s.t. \hspace{4mm} \tilde{\gamma}^{\vartriangle}=g\cdot \gamma^{\vartriangle}; \hspace{5mm} 2<i<n-1}
\end{eqnarray*}

\vspace{-3mm}
\noindent which means the congruence of $\mathrm{\gamma^{\vartriangle}}$ and $\mathrm{\tilde{\gamma}^{\vartriangle}}$.
\clearpage 

Similarly, we can prove the following theorem in terms of (6).

\textbf{Theorem 5.8, Signature-inverse Theorem.} \textit{Let $\mathrm{\gamma^{\vartriangle}, \tilde{\gamma}^{\vartriangle}\subset M^{S\Bbb{A}}}$ be two ordinary $\mathrm{\Bbb{A}}$-fine convex meshes with the same area in the one-neighborhoods of the corresponding points with zero $\mathrm{\Bbb{A}}$-curvatures. Also, let they have the same $\mathrm{\Bbb{A}}$-arc length set and 3-point $\mathrm{S\Bbb{A}}$-signatures. Then, $\mathrm{\gamma^{\vartriangle}}$ and $\mathrm{\tilde{\gamma}^{\vartriangle}}$ are congruent}.

The preceding theorem in terms of fine meshes is written as follows.

\noindent \hspace{4.3mm}\textbf{Corollary 5.9, Signature-inverse Theorem.} \textit{Let $\mathrm{\gamma^{\vartriangle}, \tilde{\gamma}^{\vartriangle}\subset M^{S\Bbb{A}}}$ denote two fine convex meshes having the same area in the one-neighborhoods of the corresponding points with zero $\mathrm{\Bbb{A}}$-curvatures. Let also they have the same $\mathrm{\Bbb{A}}$-arc length set and 3-point $\mathrm{S\Bbb{A}}$-signatures. Then, $\mathrm{\gamma^{\vartriangle}}$ and $\mathrm{\tilde{\gamma}^{\vartriangle}}$ are congruent}.

Just like the Euclidean case, more versions of Signature-inverse theorems can be presented in terms of $\mathrm{(m_{1}, m_{2})}$ $\mathrm{S\Bbb{A}}$-signatures parameterized by the identities (5), (7), and (8).
%
%
%==================================================================================================================
%
%
\vspace{-4mm}
\section{6. Conclusions}
\vspace{-2mm}

This paper has considered the new formulation for $\mathrm{\Bbb{G}}$-JINSs introduced in the first paper in this series. We first provided several counterexamples for Curvature-inverse Theorem and Signature-inverse Theorem in $\mathrm{M^{\Bbb{G}}}$, meaning that non-congruent meshes may have the same $\mathrm{\Bbb{G}}$-curvature or $\mathrm{\Bbb{G}}$-JINS. To give a numerical version of Euclidean Signature-inverse Theorem in terms of the associated JINSs, we first classified equally and unequally spaced three-point meshes with respect to their curvatures and side lengths, then, we looked for conditions that make this theorem correct. Next, we brought forward The Host Theorem to prove a simpler version of Signature-inverse Theorem for closed meshes. Finally, we went through the same process for the equiaffine case by classifying five-point ordinary meshes. 
%
%
%===========================================================================================================================
%
%
%\vspace{-5mm}
%\section{Acknowledgment} 
%\vspace{-3mm}
%
%I would like to thank Peter J. Olver for his advice and comments.
%
%
%***************************************************************************************************************************
%
%

\end{document}